\documentclass[10pt,twocolumn,twoside]{IEEEtran}

\usepackage[ansinew]{inputenc}
\usepackage{amsfonts}
\usepackage{graphics,graphicx,color,colortbl}
\usepackage{subfigure}
\usepackage{algorithm}
\usepackage{algorithmic}
\usepackage{amsmath}
\usepackage{bbold}
\usepackage{mathbbol}
\usepackage{arydshln,leftidx,mathtools}
\usepackage{bm}
\usepackage[top=72pt, bottom=54pt, left=54pt, right=54pt]{geometry}

\usepackage{amsthm}
\usepackage{breqn}
\usepackage{xr}
\usepackage{xcolor}
\usepackage{amsmath}

\newtheorem{Prop}{Proposition}[section]
\newtheorem{coll}{Corollary}[section]
\newtheorem{Defi}{Definition}[section]
\newtheorem{theor}{Theorem}[section]
\newtheorem{lem}{Lemma}[section]

\newtheorem{Rem}{Remark}[section]

\newcommand{\specnorm}[1]{{ \left\vert\kern-0.25ex\left\vert\kern-0.25ex\left\vert #1 
    \right\vert\kern-0.25ex\right\vert\kern-0.25ex\right\vert } }

\newcommand{\WideLaplacian}{\boldsymbol{\mathcal{L}}}

\begin{document}

\title{Distributed PID Control for Consensus of Homogeneous and Heterogeneous Networks}

\author{Daniel~Burbano,~\IEEEmembership{Student,~IEEE,}
        and~Mario~di Bernardo,~\IEEEmembership{Fellow,~IEEE}
\thanks{D. Burbano and M. di Bernardo are with the Department of Electrical Engineering and Information Technology, University of Naples Federico II, Naples 80125, Italy.}
\thanks{M. di Bernardo is  also with the Department of Engineering Mathematics, University of Bristol, U.K.}
}

\maketitle

\begin{abstract}
We investigate the use of distributed PID actions to achieve consensus in networks of homogeneous and heterogeneous linear systems. Convergence of the strategy is proved for both cases using appropriate state transformations and Lyapunov functions. The effectiveness of the theoretical results is illustrated via its application to a representative power grid model recently presented in the literature.
\end{abstract}

\IEEEpeerreviewmaketitle

\section{Introduction}

\IEEEPARstart The problem of driving a network of interconnected dynamical agents asymptotically towards the same state is relevant to many applications.
Examples include the design of heating, ventilation and air conditioning (HVAC) systems to obtain a constant temperature throughout a smart building \cite{Weimer2012}; distributed formation control in robotics \cite{Olfati-Saber2002I,Fax2004}; platooning of vehicles in intelligent transportation systems \cite{Coelingh2012,Naus2010}; and frequency synchronization \cite{Hill2006,Doerfler2012} in power grids and microgrids. (For a more comprehensive list of applications see \cite{Antonelli2013, OlfatFM2007} and references therein.)

From the early work reported in \cite{Saber2003}, achieving consensus in multi-agent systems and networks has become a fundamental problem in Control. The classical paradigm involves networks of simple or higher-order integrators communicating via linear diffusive coupling on an undirected network. Extensions were also presented to a number of different cases; for example those where the communication protocol is time-varying, nonlinear, or affected by switching and delays \cite{RenB2003,ZhiyuFM2005,OlfatM2004,OlfatFM2007}, the network graph is directed, and the node dynamics is nonlinear (for a review see \cite{ReCa:11,Ren2007c}).

Typically, it is assumed that the agent dynamics is either trivial or identical across the network. Thus, many of the available strategies only apply to networks of homogeneous systems in the absence of disturbances and noise. Unfortunately, in many applications this is not the case. 
Take for instance a network of power generators, as those considered in \cite{Dorfler2013}. Multiple unavoidable disturbances such as error measurements, sudden load variations, and communication failures between generators make the network highly heterogeneous. 

Some recent work addresses the problem of achieving consensus in networks with some degree of heterogeneity. For instance, in \cite{Fradkov2011}, the problem studied is of driving all the linear nodes in a homogeneous network towards a common reference trajectory (leader-follower networks) in the presence of time-varying, yet bounded, disturbances. Also, the case of heterogeneous networks has been studied in the absence of disturbances or noise both for linear \cite{Wieland2011,Kim2011,Grip2012,Wang2012} and nonlinear node dynamics \cite{ZhaoMarch,ZhongOct}. 

The pressing open problem still remains of designing strategies able to guarantee convergence of all agents towards the same solution in the presence of heterogeneity among their dynamics together with disturbances and noise. 
In this case, diffusive linear coupling is in general only able to guarantee bounded steady-state error as the coupling gain is increased \cite{Kim2012}.

To overcome some of these problems, the use of a distributed integral action to achieve consensus was proposed in the literature. For instance, PI coupling is used in \cite{CarliDZ2011} to achieve clock synchronization in networks of  discrete-time integrators. Also, distributed PI actions are exploited in \cite{Andreasson2012a} to achieve consensus in networks of simple and double integrators affected by constant disturbances.

The aim of this paper is to present a notable extension of results in \cite{CarliDZ2011,Andreasson2012a} by proposing the use of a distributed PID protocol as a simple yet effective solution to achieve consensus in networks of linear systems despite the presence of both heterogeneous node dynamics and constant disturbances.

A proof of convergence is obtained for both the  homogeneous and heterogeneous cases. Specifically, novel conditions are derived for tuning the gains of the distributed PID strategy that depend on the node dynamics and the network structure. The theoretical derivations are based on linear algebra and the use of appropriate Lyapunov functions.
The results are illustrated on the representative example of a linearized power system model which was also investigated in \cite{Simpson-Porco2013}. 

\section{Mathematical Preliminaries}
\label{sec:MP}
We denote by $\mathbf{I}_N$ the identity matrix of dimension $N\times N$; by $\mathbb{0}_{M\times N}$ a matrix of zeros of dimension $M\times N$, and by $\mathbb{1}_N$  a $N\times 1$ vector with unitary elements. The Frobenius norm of a matrix or a vector is denoted by $\left\| \cdot \right\|$ while the spectral norm of a matrix by $\specnorm{\cdot}$. A diagonal matrix, say $\mathbf{D}$, with diagonal elements $\lambda_1, \ldots,\lambda_N$ is indicated by $\mathbf{D}=\mbox{diag}\{\lambda_1,\ldots,\lambda_N\}$. The determinant of a matrix is denoted by $\det(.)$. Given two vectors ${\boldsymbol\zeta}_1$, ${\boldsymbol\zeta}_2\in \mathbb{R}^{n\times 1}$ and a matrix $\mathbf{Q}\in \mathbb{R}^{n\times n}$, from linear algebra one has \cite{Jun-Wei2011}
\begin{equation}
\label{pos:eq:1}
2\boldsymbol\zeta_1^T\mathbf{Q}^T\boldsymbol\zeta_2 \le \sigma\boldsymbol\zeta_1^T\mathbf{Q}^T\mathbf{Q}\boldsymbol\zeta_1+\frac{1}{\sigma}\boldsymbol\zeta_2^T\boldsymbol\zeta_2,\forall \sigma>0
\end{equation}
An {\em undirected graph} $\mathcal{G}$ is a pair defined by $\mathcal{G} = \left( {\mathcal{N},\mathcal{E}} \right)$ where $\mathcal{N} = \left\{ {{1},{2}, \cdots ,{N}} \right\}$ is the finite set of $N$ node indices; $\mathcal{E} \subset \mathcal{N} \times \mathcal{N}$ is the set containing the $M$ edges among nodes. Furthermore, we assume each edge has an associated weight denoted by $w_{ij} \in \mathbb{R}^{+}$. The \textit{Laplacian matrix} $\WideLaplacian(\mathcal{G})\in {\mathbb{R}^{N \times N}}$ is defined as the matrix whose elements ${\WideLaplacian_{ij}}(\mathcal{G}) = \sum\nolimits_{j = 1,j \ne i}^N {{w_{ij}}} $ if $i=j$ and $-{{w_{ij}}}$ otherwise.

\begin{lem}\cite{Zielke1983,RAH_CRJ_1987}
\label{lemm:quadratic_form}
Given a symmetric matrix $\mathbf{A}\in \mathbb{R}^{n\times n}$, then
$${\lambda _{\min }}(\mathbf{A}) {{ {{\boldsymbol\zeta}} }^T} {{\boldsymbol\zeta}} \le {{ {{\boldsymbol\zeta}} }^T}\mathbf{A} {{\boldsymbol\zeta}} \le {\lambda _{\max }}(\mathbf{A}){{ {{\boldsymbol\zeta}} }^T} {{\boldsymbol\zeta}}, \forall \boldsymbol\zeta \in \mathbb{R}^{n\times 1}$$
where ${\lambda _{\min }}(\mathbf{A})$ and  ${\lambda _{\max }}(\mathbf{A})$ denote the smallest and largest eigenvalues of $\mathbf{A}$. Moreover, $\specnorm{ \mathbf{A} } = \mathop {\max }\limits_k \left\{ {\left| {{\lambda _k}(\mathbf{A})} \right|} \right\}\leq{\left\| \mathbf{A} \right\|}$; $\lambda _k$ denoting the $k$-th eigenvalue of $\mathbf{A}$.
\end{lem}
\begin{Defi}\cite{Lu2006214}
We say that an $N\times N$ matrix $\boldsymbol{\mathcal{M}}  = [{\mathcal{M}_{ij}}], \forall i,j\in \mathcal{N}$  belongs to the set $\mathbf{\Omega}$ if it verifies the following properties:
\begin{enumerate}
	\item ${{\mathcal{M}}_{ij}} \ge 0,\,i \ne j,$ and ${{\mathcal{M}}_{ii}} =  - \sum\limits_{j = 1,j \ne i}^N {{{\mathcal{M}}_{ij}}}$,
	\item its eigenvalues in ascending order are such that $\lambda _1(\boldsymbol{\mathcal{M}})=0$ while all the others, $\lambda _k(\boldsymbol{\mathcal{M}})$, $k\in \{2,\cdots,N\}$, are real and positive.
\end{enumerate}
\label{def_laplacian}
\end{Defi}
\begin{lem}\cite{Ren2007b}
\label{lemm:simmetric_L}
Let $\mathcal{G}$ be a connected undirected graph. Then, its corresponding Laplacian matrix ${\WideLaplacian}\in \mathbf{\Omega}$ and its eigenvalues can be sorted in ascending order as $0 = {\lambda _1} < {\lambda _2} \le  \cdots  \le {\lambda_N}$.
\end{lem}
%
%
%
%
%
\subsection{Block decomposition of $\WideLaplacian$}
Next, we present a decomposition of the Laplacian matrix that will be crucial for the derivations reported in the rest of the paper. Note that such a decomposition can be avoided sometimes when proving consensus in homogeneous networks, but it is particularly useful to prove convergence in the presence of heterogeneous nodes.

As the Laplacian matrix is symmetric (the graph is undirected), according to Schur's lemma, there exists an orthogonal matrix, say $ {\mathbf{V}}:=(1/\sqrt{N}) \mathbf{U}$ such that ${\WideLaplacian}=  {\mathbf{V}} \mathbf{\Lambda}  {\mathbf{V}}^{-1}=\mathbf{U} \mathbf{\Lambda} \mathbf{U}^{-1}$, where $\mathbf{\Lambda}  = \mbox{diag}\left\{ {0,{\lambda _2}, \cdots ,{\lambda _N}} \right\}$. 
Note that the eigenvectors of ${\WideLaplacian}$ are column vectors of $\mathbf{V}$ (or equivalently row vectors of $\mathbf{V}^{-1}$).
As suggested in \cite{Kim2012}, without loss of generality, we can express the orthogonal matrix $\mathbf{U}$ and its inverse in the following block form
\begin{equation}
\label{Ublocks}
\begin{array}{l}
\mathbf{U} = \left[ {\begin{array}{*{20}{c}}
   1 & {{\mathbf{Q}_{12}}}  \\
   {{\mathbb{1}_{N - 1}}} & {{\mathbf{Q}_{22}}}  \\
\end{array}} \right],\,{\mathbf{U}^{-1}} = \left[ {\begin{array}{*{20}{c}}
   {{r_{11}}} & {{\mathbf{R}_{12}}}  \\
   {{\mathbf{R}_{21}}} & {{\mathbf{R}_{22}}}  \\
\end{array}} \right],
\end{array}
\end{equation}
where $\mathbf{Q}_{12}\in {\mathbb{R}^{1 \times (N - 1)}}$, $\mathbf{Q}_{22}\in \mathbb{R}^{{(N - 1) \times (N - 1)}}$, $r_{11}\in {\mathbb{R}}$, 
$\mathbf{R}_{12}\in {\mathbb{R}^{1 \times (N - 1)}}$, $\mathbf{R}_{21}\in {\mathbb{R}^{{(N - 1) \times 1} }}$, 
$\mathbf{R}_{22}\in {\mathbb{R}^{{(N - 1) \times (N - 1)}}}$ are blocks of appropriate dimensions and 
\begin{eqnarray}
r_{11} &=& \frac{1}{N}, \qquad \mathbf{R}_{12}= \frac{1}{N}\mathbb{1}_{N - 1}^T,\label{eq:blockdef}
\end{eqnarray}
Moreover, as ${\mathbf{V}}^{-1}={\mathbf{V}}^T$, it follows that $\mathbf{R}_{21}=(1/N)\mathbf{Q}_{12}^T$ and $\mathbf{R}_{22} = (1/N)\mathbf{Q}_{22}^T$.
Thus, we can recast $\mathbf{U}$ as 
\begin{equation}
\label{Ublocks:b}
\begin{array}{l}
\mathbf{U} =  \left[ {\begin{array}{*{20}{c}}
   1 & {N{\mathbf{R}_{21}^T}}  \\
   {{\mathbb{1}_{N - 1}}} & {N{\mathbf{R}_{22}^T}}  \\
\end{array}} \right],
\end{array}
\end{equation}
Also, since $\mathbf{U}^{-1} \mathbf{U}=\mathbf{I}_N$, the blocks in the definition of $\mathbf{U}$ and $\mathbf{U}^{-1}$ must fulfill the following conditions:
\begin{equation}
\label{prop:U:1}
{{{\mathbf{R}_{21}} + {\mathbf{R}_{22}}{\mathbb{1}_{N - 1}} = {\mathbb{0}_{\left( {N - 1} \right) \times 1}}}} 
\end{equation}
\begin{equation}
\label{prop:U:2}
{{{\mathbf{R}_{21}}{\mathbf{R}_{21}^T} + {\mathbf{R}_{22}}{\mathbf{R}_{22}^T} = \frac{1}{N}{\mathbf{I}_{N - 1}}}}
\end{equation} 
\begin{equation}
\label{prop:U:3}
{{{r_{11}}{\mathbf{R}_{21}^T} + {\mathbf{R}_{12}}{\mathbf{R}_{22}^T} = {\mathbb{0}_{1 \times \left( {N - 1} \right)}}}} 
\end{equation} 
Note that, solving (\ref{prop:U:1}) for $\mathbf{R}_{21}$ and (\ref{prop:U:3}) for $\mathbf{R}_{21}^T$, using (\ref{eq:blockdef}), we can also write 
\begin{equation}
\label{prop:U:4}
\mathbf{R}_{21}\mathbf{R}_{21}^T=\mathbf{R}_{22}\mathbb{1}_{N-1}\mathbb{1}_{N-1}^T \mathbf{R}_{22}^T
\end{equation} 
Moreover, ${\specnorm{ \mathbf{V}^{-1} } }=\sqrt{\lambda_{max}( (\mathbf{V}^{-1})^T {\mathbf{V}}^{-1})}$. Also, as $ {\mathbf{V}}^{-1}= {\mathbf{V}}^{T}$ and $\mathbf{V}\mathbf{V}^T=\mathbf{I}_{N}$ one has $\specnorm{{\mathbf{U}^{ - 1}}} = \frac{1}{{\sqrt N }}$ and therefore its block $\mathbf{R}_{22}$ is such that
\begin{equation}
\label{prop:U:6}
\specnorm{\mathbf{R}_{22}} \le {\specnorm{ \mathbf{U}^{-1} } } = {1}/{{\sqrt N }}
\end{equation}
Finally, $\left\|\mathbf{R}_{22}\mathbb{1}_{N-1}\right\|\leq \sqrt{N-1}\specnorm{\mathbf{R}_{22}}$  (see Theorem 5.6.2 of \cite{RAH_CRJ_1987}), then expressing $\mathbf{R}_{21}$ from \eqref{prop:U:1} and using (\ref{prop:U:6}), we find that 
\begin{equation}
\label{bound_R21_termsofN}
\left\| \mathbf{R}_{21}\right\|\leq\sqrt{N-1}\specnorm{\mathbf{R}_{22}}\leq\sqrt{(N-1)/N}
\end{equation}
\begin{Prop} The matrix $\mathbf{R}_{22}$ has full rank.
\label{FullRank:R22}
\end{Prop}
\begin{IEEEproof} From the definition, $\mathbf{U}$ is full rank and $\det(\mathbf{U})\ne 0$. Then, from \eqref{Ublocks:b} one has (see Prop. 2.8.3 in \cite{Bernstein2009})
$$\det(\mathbf{U})=\det(N\mathbf{R}_{22}^T-N\mathbb{1}_{N-1}\mathbf{R}_{21}^T)$$ 

Now, from \eqref{prop:U:1}, $\mathbf{R}_{21}^T=-\mathbb{1}_{N-1}^T\mathbf{R}_{22}^T$. Therefore,
 $$\det(\mathbf{U})=N\det(\mathbf{I}_{N-1}+\mathbb{1}_{N-1}\mathbb{1}_{N-1}^T)\det(\mathbf{R}_{22}^T)\ne 0$$ which implies that the rank of $\mathbf{R}_{22}$ is full.
\end{IEEEproof}
Next, we investigate some properties of the matrix $\tilde {\WideLaplacian}$ defined as $\tilde{\WideLaplacian}:= \mathbf{I}_N + \gamma \WideLaplacian$, where $\gamma \ge 0$ will be used to denote the derivative gain in the next section.
From the definition and the properties of $\WideLaplacian$, we immediately have that $\tilde {\WideLaplacian}$ is invertible and its eigenvalues are $\tilde \lambda_i = 1+\gamma \lambda_i, \forall i\in \mathcal{N}$.

\subsection{Properties of $\tilde{\WideLaplacian}$}
We can prove the following Lemma that will be useful for the convergence analysis reported in Section \ref{sec:PID}.
%
\begin{lem} 
\label{lem:LDL_eig}
If ${\WideLaplacian}\in \mathbf{\Omega}$, then $\tilde {\WideLaplacian}^{ - 1}$ has positive real eigenvalues that can be given in descending order as $1 \ge 1/(\gamma {\lambda _2} + 1) \ge  \cdots  \ge 1/(\gamma {\lambda _N} + 1)$. Moreover, the product $\tilde {{\WideLaplacian}}^{ - 1}{\WideLaplacian}$ is itself in $\mathbf{\Omega}$ and can be expressed as
\begin{equation}
\tilde {{\WideLaplacian}}^{ - 1}{\WideLaplacian}= \mathbf{U}\mathbf{\Gamma} {\mathbf{U}^{-1}},\,\,\,\ \mathbf{\Gamma}  = \left[ {\begin{array}{*{20}{c}}
0&{{\mathbb{0}_{1 \times (N - 1)}}}\\
{{\mathbb{0}_{(N - 1) \times 1}}}&{\widehat {\mathbf{\Gamma}} }
\end{array}} \right]\,
\label{eq:eig_LD}
\end{equation}
with
\begin{equation}
\label{Gamma:1}
\widehat{\mathbf{\Gamma}}  = \mbox{diag} \left\{{{\lambda _2}/\left( {\gamma{\lambda _2} + 1} \right), \cdots ,{\lambda _N}/\left( {\gamma{\lambda _N} + 1} \right)} \right\}
\end{equation} 
\end{lem}
\begin{IEEEproof} From the definition of $\tilde{{\WideLaplacian}}$ we can write $\tilde{{\WideLaplacian}} = \mathbf{U}\mathbf{U}^{-1} + \gamma\mathbf{U}\mathbf{\Lambda}\mathbf{U}^{-1}$, and letting $\mathbf{\Sigma} : = \mbox{diag}\left\{ {1,\gamma {\lambda _2} + 1, \cdots ,\gamma {\lambda _N} + 1} \right\}$ yields $\tilde{{\WideLaplacian}} = \mathbf{U}\mathbf{\Sigma}\mathbf{U}^{-1}$. Note that this is the eigen-decomposition of a symmetric matrix and $\tilde{{\WideLaplacian}}^{-1} = \mathbf{U}\mathbf{\Sigma}^{-1}\mathbf{U}^{-1}$ where
\begin{equation}
\label{Sigma_inv}
{\mathbf{\Sigma}}^{-1}:=\mbox{diag}\{ {1,{1}/\left( {\gamma{\lambda _2} + 1} \right), \cdots ,{1}/\left( {\gamma{\lambda _N} + 1} \right)} \}
\end{equation}
Thus, we obtain $\tilde{{\WideLaplacian}}^{-1}{\WideLaplacian}=\mathbf{U}{\mathbf{\Gamma}}\mathbf{U}^{-1}$ with ${\mathbf{\Gamma}}:=\mathbf{\Sigma}\mathbf{\Lambda}$ and the proof is complete. 
\end{IEEEproof}
Also, rewriting $\tilde {\WideLaplacian}^{-1}$ in block form, we have
\begin{equation}
\label{inv_lapla}
{{\tilde {{\WideLaplacian}}}^{ - 1}} = \left[ {\begin{array}{*{20}{l}}
{\widehat {l}_{11}}&{\widehat {\WideLaplacian}_{12}}\\
{\widehat {\WideLaplacian}_{21}}&{\widehat {\WideLaplacian}_{22}}
\end{array}} \right]
\end{equation}
where ${\widehat {l}_{11}}\in\mathbb{R}$, ${\widehat {\WideLaplacian}_{12}} \in {\mathbb{R}^{1 \times N - 1}}$, ${\widehat {\WideLaplacian}_{21}} \in {\mathbb{R}^{N - 1 \times 1}}$, and ${\widehat {\WideLaplacian}_{22}} \in {\mathbb{R}^{N - 1 \times N - 1}}$ are blocks of appropriate dimension. Moreover, from Def. \ref{def_laplacian}, $\WideLaplacian{\mathbb{1}_N} = {\mathbb{0}_{N\times 1}}$; therefore, $\tilde {\WideLaplacian}{\mathbb{1}_N} = {\mathbb{1}_N}$. Hence, multiplying both sides by ${\tilde {\WideLaplacian}^{ - 1}}$ we have ${\tilde {\WideLaplacian}^{ - 1}}\tilde{\WideLaplacian}{\mathbb{1}_N} =  {\tilde {\WideLaplacian}^{ - 1}}{\mathbb{1}_N}$ so that ${\tilde {\WideLaplacian}^{ - 1}}{\mathbb{1}_N} = {\mathbb{1}_N}$. Therefore, the blocks in \eqref{inv_lapla} must satisfy the conditions:
\begin{eqnarray}
\label{prop:L:1}
\widehat {\WideLaplacian}_{12}{\mathbb{1}_{N - 1}} &=& { {\widehat {\WideLaplacian}_{21}^T} }{\mathbb{1}_{N - 1}} =1 - \widehat {l}_{11}\\
\label{prop:L:2}
\widehat {\WideLaplacian}_{22}{\mathbb{1}_{N - 1}} &=& {\mathbb{1}_{N - 1}} - \widehat {\WideLaplacian}_{21}\\
\label{prop:L:3}
{\mathbb{1}_{N - 1}}\mathbb{1}_{N - 1}^T\widehat {\WideLaplacian}_{22} &=& {\mathbb{1}_{N-1}}\mathbb{1}_{N - 1}^T - {\mathbb{1}_{N-1}}{\widehat {\WideLaplacian}_{12}}\\
\label{prop:L:4}
\widehat{ {{\mathbf{H}}}}{\mathbb{1}_{N - 1}} &=&  {{{\widehat l}}_{11}{\mathbb{1}_{N - 1}} - {{\widehat {\WideLaplacian}}}_{21}}
\end{eqnarray}
where
\begin{equation}
\label{eq:hat_H}
\widehat{\mathbf{H}} :=  {\widehat {\WideLaplacian}_{22} - {\mathbb{1}_{N - 1}}\widehat{\WideLaplacian}_{12}}
\end{equation}
From Lemma \ref{lem:LDL_eig} we have that ${\mathbf{U}^{-1}}\tilde {\WideLaplacian}^{ - 1}\mathbf{U}= {\mathbf{\Sigma}} ^{ - 1}$, then applying some block operations and using properties (\ref{prop:U:1})-(\ref{prop:U:3}) and (\ref{prop:L:1})-(\ref{prop:L:3}) one also has (see Appendix \ref{Appendix_I})
\begin{equation}
\label{prop:L:5}
\begin{array}{l}
\frac{1}{N}  {\mathbf{\widehat\Sigma}^{-1}} ={\mathbf{R}_{22}} \left( \widehat{ {{\mathbf{H}}}} + ({{\widehat l}_{11}}{\mathbb{1}_{N - 1}} - \widehat {\WideLaplacian}_{21})\mathbb{1}_{N - 1}^T\right)\mathbf{R}_{22}^T
\end{array}
\end{equation} 
where 
\begin{equation}
\label{bar_Sigma}
{\mathbf{\widehat \Sigma}^{-1}} := \mbox{diag}\left\{ {1/(\gamma {\lambda _2} + 1), \cdots ,1/(\gamma {\lambda _N} + 1)} \right\}
\end{equation} 
%
%
\section{Problem Formulation}
\label{Sec:III}
We consider a group of $N$ nodes governed by heterogeneous first-order linear dynamics of the form
\begin{equation}
\label{eq:sys:1}
{{\dot x}_i}(t) = {\rho _i}{x_i}(t) + {\delta _i} + {u_i}(t),\,\,\,i \in \mathcal{N}	
\end{equation}
where $x_i(t)\in \mathbb{R}$ represents the state of the $i$-th agent, ${\rho _i}\in \mathbb{R}$ is the agent pole determining its uncoupled dynamics,  ${\delta _i}\in \mathbb{R}$ is some constant disturbance (or constant external input) acting on each node, and $u_i(t)\in \mathbb{R}$ is the distributed control input through which agent $i$ communicates with its neighboring agents. Note that without any control input, the node dynamics can either be stable (${\rho _i} < 0$) or unstable (${\rho _i} > 0$), while the constant term $\delta _i$ can be used to represent different quantities in applications. For example, it can model constant power injections in power grids \cite{Doerfler2012} or noise in minimal models of flocks of birds \cite{Young2013}. 

Let $\mathbf{x}(t) = \left[ {{x_1}(t), \ldots, {x_N}(t)} \right]^T$ be the stack vector of all node states and define $\mathfrak{C}$ as the {\em consensus manifold} 
$$\mathfrak{C} := \left\{ {\left. {\mathbf{x} \in \mathbb{R}^N\,\,} \right|\,\left|{x_j}(t) - {x_i}(t)\right| = 0\,}, \forall i,j\in \mathcal{N}, i\neq j \right\}$$ 
\begin{Defi} (\textit{Admissible consensus}) 
\label{def:1}
The network of $N$ heterogeneous agents described by (\ref{eq:sys:1}) is said to reach admissible consensus if, for any set of initial conditions $x_i(0)=x_{i0}$,
\[\mathop {\lim }\limits_{t \to \infty } {\mathbf{x}(t)}\in\mathfrak{C},\quad \left|{u_i(t)}\right|  <  + \infty, \quad \forall t\geq0,\ i\in \mathcal{N}\]
If, instead ${\lim _{t \to \infty }}\left| {{x}_j}(t) - {{x}_i}(t) \right| \le \varepsilon$ for $\varepsilon>0$ and $\left|{u_i(t)}\right|  <  + \infty, \quad \forall t\geq0,\ \forall i,j\in \mathcal{N}$, the network is said to achieve $\varepsilon$-admissible consensus.
\end{Defi}
In this paper we study the case where, rather than communicating via a classical proportional (diffusive) coupling, agents in the network are coupled via the PID consensus protocol given by
\begin{equation}
\label{eq:cont:2}
{u_i}(t) =  - \sum\limits_{j = 1}^N {{\WideLaplacian_{ij}}} \left( {\alpha {x_j}(t) + \beta \int\limits_0^t {{x_j}(\tau)d\tau}  + \gamma {{\dot x}_j}(t)} \right)
\end{equation}
where ${\WideLaplacian}_{ij}$ are the elements of the network Laplacian and $\alpha > 0, \beta \geq 0$, $\gamma \geq 0$ are the gains determining the strength of the proportional, integral and derivative actions, respectively. 
Defining $\mathbf{P}:= \mbox{diag}\left\{ {{\rho _1}, \cdots ,{\rho _N}} \right\}$,  $\mathbf{\Delta} := {\left[ {{\delta _1}, \cdots ,{\delta _N}} \right]^T}$, and the stack vector of integral states 
\begin{equation}
\label{integral:term}
{\mathbf{z}}(t) =\left[ {z_1}(t), \ldots ,{z_N}(t) \right]^T:= - \beta \tilde {\WideLaplacian}^{-1}\WideLaplacian\int_0^t {\mathbf{x}(\tau )d\tau }
\end{equation}
the overall dynamics of the closed-loop network can be written as (using the notation introduced in Sec. \ref{sec:MP})
\noindent
\begin{equation}
\label{eq:PID:1}
\left[ {\begin{array}{*{20}{c}}
{\dot {\mathbf{x}}(t)}\\
{\dot {\mathbf{z}}(t)}
\end{array}} \right] = \underbrace {\left[ {\begin{array}{*{20}{c}}
{\mathbf{A}_{1}}&{{\mathbf{I}_N}}\\
{\mathbf{A}_{2}} & \mathbb{0}
\end{array}} \right]}_{\mathbf{A}}\left[ {\begin{array}{*{20}{c}}
{\mathbf{x}(t)}\\
{\mathbf{z}(t)}
\end{array}} \right] + \left[ {\begin{array}{*{20}{c}}
{\tilde {\WideLaplacian}^{-1}\mathbf{\Delta} }\\
\mathbb{0}
\end{array}} \right]
\end{equation}
where $\mathbf{A}_{1}:=\tilde {\WideLaplacian}^{-1}\left( {\mathbf{P} - \alpha \WideLaplacian} \right)$ and $\mathbf{A}_{2}:=- \beta \tilde {\WideLaplacian}^{-1}\WideLaplacian$. Then, the problem is finding conditions on the control gains $\alpha$, $\beta$ and $\gamma$, the network structure and node dynamics such that the closed-loop network \eqref{eq:PID:1} achieves admissible consensus.
\begin{Prop}
\label{equil_point}
Network \eqref{eq:PID:1} has a unique equilibrium given by ${\mathbf{x}^*}:=x_{\infty}\mathbb{1}_{N}$ with  $x_{\infty}:=-{{\sum\nolimits_{k = 1}^N {{\delta _k}} }} / {{\sum\nolimits_{k = 1}^N {{\rho _k}} }}$, and ${\mathbf{z}^* }:=- \tilde {\WideLaplacian}^{-1}( {\mathbf{P}{\mathbf{x}^* } + \mathbf{\Delta} } )$.
\end{Prop}
\begin{IEEEproof}
The proof follows immediately by setting the left-hand side of \eqref{eq:PID:1} to zero  and noticing that  $\tilde {\WideLaplacian}^{ - 1}\WideLaplacian \in \mathbf{\Omega}$ (see Lemma \ref{lem:LDL_eig}) so that $\mathbf{x}^* = a\mathbb{1}_N$, $\forall a \in \mathbb{R}$ and ${\mathbf{z}^* } =  - \tilde {\WideLaplacian}^{-1}\left( {a\mathbf{P}\mathbb{1}_N + \mathbf{\Delta} } \right)$. By definition \eqref{integral:term}, and from the fact that $\tilde {\WideLaplacian}^{-1}\WideLaplacian\in{\mathbf{\Omega}}$, we also have $\mathbb{1}_N^T\mathbf{z}(t) = 0$, then $\mathbb{1}_N^T\mathbf{z}^*  = 0$ and we obtain 
$$a=  - \mathbb{1}_N^T\mathbf{\Delta} /\mathbb{1}_N^T\mathbf{P}{\mathbb{1}_N} =  - \left(\sum\nolimits_{k = 1}^N {{\delta _k}} \right) \left(\sum\nolimits_{k = 1}^N {{\rho _k}}\right)^{-1}$$ 
Hence, setting $x_{\infty}=a$ completes the proof.
\end{IEEEproof}
Note that, results presented in this paper significantly extend previous ones in the literature  where all nodes are simple integrators, $\mathbf{P}=\mathbb{0}_{N\times N}$, \cite{CarliDZ2011} or share identical stable dynamics, $\mathbf{P}=p\mathbf{I}_N$, $p<0$ \cite{Andreasson2012a}. Specifically, we prove that there exists an $\alpha$ such that $\mathbf{A}_{1}$ is Hurwitz if $\mathbf{P}$ is a generic diagonal matrix with negative trace, i.e. even if nodes have heterogeneous dynamics and some are possibly unstable.
%
%
%
%
%

\section{Convergence analysis}
\label{sec:PID}
To prove convergence of the closed-loop network \eqref{eq:PID:1}, we choose to analyse the transverse stability of the consensus manifold $\mathfrak{C}$. We split the proof into two stages. Firstly, the system describing the transverse dynamics to the consensus manifold is derived and some of its generic properties are described. Secondly, the two cases of homogeneous and heterogeneous nodes are treated separately  in order to complete the proof of convergence. It is shown that the distributed integral action can be effectively used to reject constant disturbances, while the distributed derivative action decreases bounds on the integral terms.
\subsection{Step 1: Transverse Dynamics}
To study convergence to the consensus manifold, we consider the state transformation $\mathbf{x}^\perp (t)  = {\mathbf{U}^{-1}}\mathbf{x}(t)$. Indeed, using the block representation of $\mathbf{U}^{-1}$ in (\ref{Ublocks}) and letting ${\bar {\mathbf{x}}}^\perp (t) = \left[ x_2^\perp (t), \ldots , x_N^\perp (t) \right]$, ${\bar{\mathbf{x}}}(t) = \left[x_2(t), \ldots ,x_N(t) \right]$ we obtain
\begin{subequations}
\label{eq:transformation_U}
\begin{alignat}{3}
    \label{eq:transformation_Ua}
    {x_1^\perp(t)} &= {r_{11}}{x_1}(t) + {\mathbf{R}_{12}} {\bar{\mathbf{x}}}(t) \\
    \label{eq:transformation_Ub}
    {\bar{\mathbf{x}}^\perp(t)} &= \mathbf{R}_{21} x_1(t) + \mathbf{R}_{22} \bar{\mathbf{x}}(t)
\end{alignat}
\end{subequations}
Note that by adding and subtracting the term $\mathbf{R}_{22}x_1(t)\mathbb{1}_{N-1}$ to (\ref{eq:transformation_Ub}), and using property (\ref{prop:U:1}), one has $${{\bar {\mathbf{x}}^\perp}} (t) = {\mathbf{R}_{22}}\left( {\bar{\mathbf{x}}(t) - {x_1}(t)\mathbb{1}_{N-1}} \right)$$ It is important to highlight that ${{\bar {\mathbf{x}}^\perp}} (t) =\mathbb{0}$ if and only if $ {\bar{\mathbf{x}}(t) - {x_1}(t)\mathbb{1}_{N-1}} =\mathbb{0}$ since ${\mathbf{R}_{22}}$ is full rank (Proposition \ref{FullRank:R22}). Then, admissible consensus is achieved if $\lim_{t \to \infty }{{\bar {\mathbf{x}}^\perp}} (t)=\mathbb{0}$ and $\left\| {{\mathbf{z}}(t)} \right\| < +\infty, \forall t>0$. 

Now, recasting \eqref{eq:PID:1}, in the new coordinates $\mathbf{x}^\perp (t)$ and $\mathbf{z}^\perp (t)  = {\mathbf{U}^{-1}}\mathbf{z}(t)$, and using Lemma \ref{lem:LDL_eig}, we get
\begin{subequations}
\label{eq:cont:PID_Proj}
\begin{alignat}{3}
\label{eq:cont:PID_Proj:a}
\dot{\mathbf{x}}^\perp(t) &= \left( {{\mathbf{\Psi}}  - \left[ {\begin{array}{*{20}{c}}
0&\mathbb{0}\\
\mathbb{0}&{ \alpha  {\widehat{\mathbf{\Gamma}} }}
\end{array}} \right]} \right){\mathbf{x}}^\perp(t) + {\mathbf{z}}^\perp(t) + \tilde{\mathbf{\Delta}}  \\
\label{eq:cont:PID_Proj:b}
\dot{\bar{\mathbf{z}}}^\perp(t)  &=  - \beta {\widehat{\mathbf{\Gamma}} }{{ {\mathbf{\bar x}}^\perp}}(t) 
\end{alignat}
\end{subequations}
where ${ {{\bar {\mathbf{z}}}}}^\perp (t) = \left[z_2^\perp (t), \ldots , z_N^\perp (t) \right]$ and $\mathbf{\Psi}  := {\mathbf{U}^{ - 1}}\tilde {\WideLaplacian}^{ - 1}\mathbf{P}\mathbf{U}$. (Note that the equation for ${z}_1^\perp(t)$ can be neglected as it has trivial dynamics with null initial conditions and represents an uncontrollable and unobservable state.) Furthermore, $\tilde{\mathbf{\Delta}} := \mathbf{U}^{-1}\tilde {\WideLaplacian}^{ - 1}\mathbf{\Delta}=[\mathbf{q}^T\,\,\,\widehat{\mathbf{R}}^T]^T\mathbf{\Delta}$, where $\mathbf{q}\in \mathbb{R}^{1 \times N}$, $\widehat{\mathbf{R}}\in \mathbb{R}^{(N-1) \times N}$ are given by
\[\left[ \begin{array}{l}
\mathbf{q}
\\
\widehat{\mathbf{R}}
\end{array} \right]=\left[ {\begin{array}{*{20}{c}}
{r_{11}\widehat l_{11} + \mathbf{R}_{12}\widehat{\WideLaplacian}_{21}}&{r_{11}\widehat{\WideLaplacian}_{12}+\mathbf{R}_{12}\widehat{\WideLaplacian}_{22}}\\
{\mathbf{R}_{21}\widehat l_{11} + \mathbf{R}_{22}\widehat{\WideLaplacian}_{21}}&{\mathbf{R}_{21}\widehat{\WideLaplacian}_{12} + \mathbf{R}_{22}\widehat{\WideLaplacian}_{22}}
\end{array}} \right]
\]
From the definition of $\mathbf{R}_{12}$ in \eqref{eq:blockdef} and using \eqref{prop:L:1} we obtain 
\begin{equation}
\mathbf{q}=[r_{11}\,\, \mathbf{R}_{12}]=\frac{1}{N}\mathbb{1}_{N}^T
\end{equation}
Also, it follows from \eqref{prop:U:1} that $\mathbf{R}_{21}=-\mathbf{R}_{22}\mathbb{1}_{N-1}$ and again using \eqref{prop:L:1} one has
$\widehat {\mathbf{R}}=[\mathbf{R}_{22}(\widehat{\WideLaplacian}_{21}-\widehat{l}_{11}\mathbb{1}_{N-1})\,\,\,\,\, \mathbf{R}_{22}(\widehat{\WideLaplacian}_{22}-\mathbb{1}_{N-1}\widehat{\WideLaplacian}_{12})]$. Thus, using \eqref{prop:L:4} we obtain
\begin{equation}
\label{eq:barR}
\widehat {\mathbf{R}} = {\mathbf{R}_{22}}\widehat{\mathbf{H}} [ - {\mathbb{1}_{N -1}}\,\,\,\,\, \mathbf{I}_{N-1}]
\end{equation}
Moreover, the matrix $\mathbf{\Psi}$ is a block matrix that can be expressed as
\begin{equation}
\label{eq:cont:Psieq}
\begin{split}
\mathbf{\Psi}  &= {\mathbf{U}^{ - 1}}\tilde {\WideLaplacian}^{ - 1}\mathbf{P}\mathbf{U}=\left[ {\begin{array}{*{20}{c}}
{{\psi _{11}}}&{{\mathbf{\Psi} _{12}}}\\
{{\mathbf{\Psi} _{21}}}&{{\mathbf{\Psi} _{22}}}
\end{array}} \right]\\
&={\mathbf{U}^{ - 1}}\left[ {\begin{array}{*{20}{l}}
{\widehat {l}_{11}}&{\widehat {\WideLaplacian}_{12}}\\
{\widehat {\WideLaplacian}_{21}}&{\widehat {\WideLaplacian}_{22}}
\end{array}} \right]\left[ {\begin{array}{*{20}{c}}
{{\rho _1}}&{{\mathbb{0}_{1 \times (N - 1)}}}\\
{{\mathbb{0}_{(N - 1) \times 1}}}&{\widehat {\mathbf{P}}}
\end{array}} \right]\mathbf{U}
\end{split}
\end{equation}
with $\widehat{\mathbf{P}} := \mbox{diag}\left\{ {{\rho_2}, \cdots ,{\rho_N}} \right\}$.

Using properties (\ref{prop:U:1})-(\ref{prop:U:4}) and (\ref{prop:L:1})-(\ref{prop:L:3}), some algebraic manipulation yields [see Appendix \ref{Appendix_II} for the derivation]
\begin{eqnarray}
\label{eq:PSI:11}
\psi _{11} &:=& (1/N)\sum\nolimits_{k = 1}^N {{\rho _k}}\\
\label{eq:PSI:12}
{\mathbf{\Psi}_{12}} &:=& \bar {\boldsymbol\rho} {\mathbf{R}_{22}^T}\\
\label{eq:PSI:21}
{{ {\mathbf{ \Psi}} }_{21}} &:=& {\mathbf{R}_{22}} \widehat{\mathbf{H}} {{\bar {\boldsymbol\rho} }^T}\\
\label{eq:PSI:22}
{{ {\mathbf{ \Psi}} }_{22}} &:=& N{\mathbf{R}_{22}} \widehat{\mathbf{H}} \left( {\mathbf{\bar P} + {\rho _1}{\mathbb{1}_{N - 1}}\mathbb{1}_{N - 1}^T} \right)\mathbf{R}_{22}^T
\end{eqnarray}
 where 
\begin{eqnarray}
\label{eq:cont:ro}
\bar {{\boldsymbol\rho}} &:=& \left[ {{\rho _2} - {\rho _1}, \cdots ,{\rho _N} - {\rho _1}} \right]
\end{eqnarray}
%
Finally, shifting the origin of \eqref{eq:cont:PID_Proj} via the further state transformation
\begin{equation}  
\label{eq:cont:shif_PID}
\left[ {\begin{array}{*{20}{c}}
{{{\hat x }_1}(t)}\\
{\hat {\mathbf{x}} (t)}\\
{\hat {\mathbf{z}}(t)}
\end{array}} \right] = \left[ {\begin{array}{*{20}{c}}
{{{ x }_1^\perp}(t)}\\
{{{{\bar {\mathbf{x}}^\perp}} (t)}}\\
{{{{\bar {\mathbf{z}}^\perp}} (t)}}
\end{array}} \right] + \left[ {\begin{array}{*{20}{c}}
{(1/{\psi _{11}})\mathbf{q}^T}\\
{{\mathbb{0}_{N - 1 \times N}}}\\
{\widehat {\mathbf{R}} - {(1/\psi _{11}){{{{\mathbf{\Psi}} _{21}}}}{{{}}}}\mathbf{q}^T}
\end{array}} \right]\mathbf{\Delta} 
\end{equation}
we obtain
\begin{equation}
 \label{eq:final:PID}
 \left[ {\begin{array}{*{20}{c}}
{{{\dot {\hat x} }_1}(t)}\\
{\dot {\hat {\mathbf{x}}} (t) }\\
{\dot {\hat {\mathbf{z}}}(t)}
\end{array}} \right] = \left[ {\begin{array}{*{20}{c}}
{\psi _{11}}&{\mathbf{\Psi}_{12}}&{\mathbb{0}}\\
{{\mathbf{\Psi}}_{21}}&{{{{\mathbf{\Psi}} _{22}} - \alpha \widehat {\mathbf{\Gamma}} }}&{\mathbf{I}_{N-1}}\\
{0}&{- \beta \widehat {\mathbf{\Gamma}} }&{\mathbb{0}}
\end{array}} \right]\left[ {\begin{array}{*{20}{c}}
{{{ {\hat x} }_1}(t)}\\
{ {\hat {\mathbf{x}}} (t) }\\
{ {\hat {\mathbf{z}}} (t) }
\end{array}} \right]
\end{equation}
Now, we can address the admissible consensus problem for \eqref{eq:PID:1} in terms of finding conditions on $\alpha$, $\beta$ and $\gamma$ that render the origin a stable equilibrium point of  \eqref{eq:final:PID}. 

Before studying this problem in the two cases of homogeneous or heterogeneous nodes, we first obtain an upper bound on the integral states that apply to both. 
By definition we have that ${\mathbf{z}}(t) = \mathbf{U}{ \mathbf{z}}^\perp(t)$, thus
\begin{subequations}
\label{for_bound_general}
\begin{alignat}{3}
{z_1}(t) &= {z_1^\perp}(t) + N\mathbf{R}_{21}^T{\bar {\mathbf{z}}}^\perp(t)\\
{\bar {\mathbf{z}}}(t) &= {\mathbb{1}_{N - 1}}{z_1^\perp}(t) + N\mathbf{R}_{22}^T{\bar {\mathbf{z}}}^\perp(t)
\end{alignat}
\end{subequations}
where ${{{\bar {\mathbf{z}}}}}(t) = \left[z_2(t), \ldots ,z_N(t) \right]$. Neglecting $z_1^\perp(t)$ for the same reason given above and using \eqref{bound_R21_termsofN}, we find from \eqref{for_bound_general} that $\left\| {{z_1}(t)} \right\| \le \sqrt{N(N-1)} \left\| {{{{\bar {\mathbf{z}}}}^ \bot(t) }} \right\|$ and $\left\| {{\bar {\mathbf{z}}}(t)} \right\| \le \sqrt{N}\left\| {{{{{\bar {\mathbf{z}}}}}^ \bot(t) }} \right\|$. Hence, we can conclude that the integral action remains bounded and $\left\| {{\mathbf{z}}(t)} \right\| \le \sqrt{N(N-1)}\left\| {{{{\bar {\mathbf{z}}}}^ \bot(t) }} \right\|$. 
Therefore, asymptotically, we have
\begin{equation}
\label{bound_z:a}
{z_\infty}:=\mathop {\lim }\limits_{t \to \infty } {\left\| {{\mathbf{z}}(t)} \right\|} \le \sqrt{N(N-1)}\mathop {\lim }\limits_{t \to \infty } {\left\| {{{{\bar {\mathbf{z}}}}^ \bot(t) }} \right\|}
\end{equation}
An upper bound for $z_\infty$ can be obtained by noticing that, if the origin of \eqref{eq:final:PID} is stable, then 
\begin{equation}
 \label{convergence:point_xi}
 \mathop {\lim }\limits_{t \to \infty }{{{\bar {\mathbf{z}}}}^ \bot(t) }= (1/{\psi _{11}}{{ {{{\mathbf{\Psi}} _{21}}}  }\mathbf{q}^T - \widehat {\mathbf{R}} })\mathbf{\Delta} 
 \end{equation}
so that
\begin{equation}
\label{bound_z:b}
\mathop {\lim }\limits_{t \to \infty } {\left\| {{{{\bar {\mathbf{z}}}}^ \bot(t) }} \right\|}\le ( 1/{\left|\psi _{11}\right|} \left\| {\mathbf{\Psi}} _{21} \right\| \left\| \mathbf{q} \right\| + {\specnorm{   \widehat{\mathbf{R}} }}) \left\| \mathbf{\Delta} \right\|
\end{equation}
Then, using \eqref{eq:barR} and \eqref{prop:U:6} yields
\begin{dmath}
\label{norm_R:a}
{\specnorm{ \widehat{\mathbf{R}} }}    \le    \specnorm{ {{\mathbf{R}_{22}}} }{\specnorm {\widehat{\mathbf{H}}}}\specnorm { {[ - {\mathbb{1}_{N - 1}}\,\,{\mathbf{I}_{N - 1}}]} }    \le   \sqrt {N} \specnorm{ {{\mathbf{R}_{22}}} }{\specnorm {\widehat{\mathbf{H}}}}
\le {\specnorm {\widehat{\mathbf{H}}}}
\end{dmath}
Also, from \eqref{eq:PSI:21} we have that
$\left\|  {\mathbf{\Psi} }_{21} \right\|\leq \specnorm{\mathbf{R}_{22}}\left\|\bar {\boldsymbol\rho}\right\|{\specnorm {\widehat{\mathbf{H}}}}$. Then, using property \eqref{prop:U:6} again, from \eqref{bound_z:b}, we can write
\begin{equation}
\label{bound_z:c}
\begin{array}{l}
\mathop {\lim }\limits_{t \to \infty } {\left\| {{{{\bar {\mathbf{z}}}}^ \bot(t) }} \right\|} \le {\specnorm {\widehat{\mathbf{H}}}} \left( 1 + \frac{\left\|\bar {\boldsymbol\rho}\right\|}{N\left| {{\psi _{11}}} \right|}  \right)\left\| \mathbf{\Delta}  \right\|
\end{array}
\end{equation}
Finally, combining \eqref{bound_z:c} with \eqref{bound_z:a} yields
\begin{equation}
\label{bound_z}
{z_\infty} \le \sqrt{N(N-1)}{\specnorm {\widehat{\mathbf{H}}}}\left( 1  + \frac{\left\|\bar {\boldsymbol\rho}\right\|}{N\left| \psi _{11} \right|}\right) \left\|\mathbf{\Delta}\right\|
\end{equation}
Note that ${\specnorm {\widehat{\mathbf{H}}}}$ is a function of $\gamma$ and therefore varying $\gamma$ controls the upper bound on $\mathbf{z}(t)$ and  consequently can be used to reduce the control effort. Specifically, we can prove the following result.
\begin{Prop} The spectral norm of $\widehat{ {\mathbf{H}} }$ can be upper bounded as
\begin{equation}
\label{bound_H}
{\specnorm {\widehat{\mathbf{H}}} } \le \frac{N}{\gamma \lambda_2+1}
\end{equation}
\label{prop:Hmatr}
\end{Prop}
\begin{IEEEproof} From Lemma \ref{lem:LDL_eig} we have that $\tilde {\WideLaplacian}^{-1} = \mathbf{U}{\mathbf{\Sigma} ^{ - 1}} \mathbf{U}^{ - 1}$. Then, using for each matrix its block representation as shown in Appendix \ref{Appendix_I}, we have
\begin{eqnarray}
\label{til_L:11}
\widehat {l}_{11} &=& r_{11}+N\mathbf{R}_{21}^T\widehat{\mathbf{\Sigma}}^{-1}\mathbf{R}_{21}\\
\label{til_L:12}
\widehat {\WideLaplacian}_{12} &=& \mathbf{R}_{12} + N\mathbf{R}_{21}^T\widehat{\mathbf{\Sigma}}^{-1}\mathbf{R}_{22}\\
\label{til_L:21}
\widehat {\WideLaplacian}_{21} &=& r_{11}\mathbb{1}_{N-1}+N\mathbf{R}_{22}^T\widehat{\mathbf{\Sigma}}^{-1}\mathbf{R}_{21}\\
\label{til_L:22}
\widehat {\WideLaplacian}_{22} &=& \mathbb{1}_{N-1}\mathbf{R}_{12}+N\mathbf{R}_{22}^T\widehat{\mathbf{\Sigma}}^{-1}\mathbf{R}_{22}
\end{eqnarray}
Replacing \eqref{til_L:12}, \eqref{til_L:22} in \eqref{eq:hat_H} and taking into account that $\mathbf{R}_{12}=1/N\mathbb{1}_{N - 1}^T$ and $\mathbf{R}_{21}^T=-\mathbb{1}_{N - 1}^T\mathbf{R}_{22}^T$ [from \eqref{prop:U:1}], we have
\begin{dmath}
\label{norm_H:a}
\widehat{\mathbf{H}} = N\mathbf{R}_{22}^T\widehat{\mathbf{\Sigma}}^{-1}\mathbf{R}_{22} - N\mathbb{1}_{N - 1}\mathbf{R}_{21}^T\widehat{\mathbf{\Sigma}}^{-1}\mathbf{R}_{22} 
= N\mathbf{R}_{22}^T\widehat{\mathbf{\Sigma}}^{-1}\mathbf{R}_{22}  + N\mathbb{1}_{N - 1}\mathbb{1}_{N - 1}^T\mathbf{R}_{22}^T\widehat{\mathbf{\Sigma}}^{-1}\mathbf{R}_{22} 
= N(\mathbf{I}_{N-1}+\mathbb{1}_{N - 1}\mathbb{1}_{N - 1}^T)\mathbf{R}_{22}^T\widehat{\mathbf{\Sigma}}^{-1}\mathbf{R}_{22}  
\end{dmath}
From \eqref{bar_Sigma} we have that $\specnorm {\widehat{\mathbf{\Sigma}}^{-1}}=1/(\gamma\lambda_2+1)$ and ${\specnorm {\widehat{\mathbf{H}}}} \le N^2\specnorm {\mathbf{R}_{22}}^2(1/(\gamma\lambda_2+1))$. Therefore, using \eqref{prop:U:6} we obtain \eqref{bound_H} and the proof is complete.
\end{IEEEproof} 
%
%
\subsection{Step 2a: Homogeneous Node Dynamics}
We first complete the proof of convergence for the homogeneous case, that is we assume all nodes share identical uncoupled dynamics.
%
%
\begin{theor}
\label{TH:PID_integrators}
The closed-loop network \eqref{eq:PID:1} with $\rho_i=-\rho^*$, $\rho^*\in\mathbb{R}^+$ $\forall i \in \mathcal{N}$ achieves admissible consensus for any positive value of $\alpha$, $\beta$ and $\gamma$. Moreover, all node states converge asymptotically to ${x_\infty} = (1/N)\sum\nolimits_{k = 1}^N {{\delta _k}/} {\rho ^*}$ with 
\begin{equation}
\label{boundPI_heter}
{z}_\infty \le  \frac{\sqrt{N^3(N-1)}}{\gamma\lambda_2+1}\left\| \mathbf{\Delta}  \right\|
\end{equation}
\end{theor}

\begin{IEEEproof} 
Firstly, note that when all nodes share the same dynamics we have  $\bar {{\boldsymbol\rho} }=\mathbb{0}_{1\times (N-1)}$ in \eqref{eq:cont:ro}. Consequently, $\mathbf{\Psi}_{12}$ and $\mathbf{\Psi}_{21}$ as defined by
\eqref{eq:PSI:12} and \eqref{eq:PSI:21} are both null vectors so that the dynamics of $\hat x_1$ in (\ref{eq:final:PID}) is independent from all the other variables, and converges to zero.

We can then study independently, the dynamics of ${\hat {\mathbf{x}}} (t)$ and ${\hat {\mathbf{z}}} (t)$ by considering the transverse candidate Lyapunov function
\begin{equation}
\label{eq:lyap}
V\left(\hat{\mathbf{x}},\hat{{z}}\right) = (\beta/2){\hat {\mathbf{x}} ^T}\widehat {\mathbf{\Gamma}} \hat {\mathbf{x}}  + \left( {1/2} \right){\hat {{\mathbf{z}}}^T}\hat {{\mathbf{z}}}
\end{equation}
which is positive definite and radially unbounded for any $\beta,\gamma>0$. Then differentiating $V$ along the trajectories of (\ref{eq:final:PID}) one has $\dot {V} =  - \beta{{\hat {\mathbf{x}} }^T}\widehat {\mathbf{\Gamma}} \left( {\alpha \widehat {\mathbf{\Gamma}}  - {{ {\mathbf{\Psi}} }_{22}}} \right)\hat {\mathbf{x}}$. As all poles are identical, $\widehat{\mathbf{P}}=-\rho^*\mathbf{I}_{N-1}$; hence, expression \eqref{eq:PSI:22} can be written as
%
%
\begin{equation}
\label{Psi:22:eq}
\begin{array}{l}
{{ {\mathbf{\Psi}}}_{22}} =  {-\rho ^*}N ({\mathbf{R}_{22}}\widehat{ {{\mathbf{H}}}}\mathbf{R}_{22}^T  +{\mathbf{R}_{22}}\widehat{ {{\mathbf{H}}}}{\mathbb{1}_{N - 1}}{\mathbb{1}_{N - 1}^T}\mathbf{R}_{22}^T)
\end{array}
\end{equation}
and, using properties \eqref{prop:L:4} and \eqref{prop:L:5}, we obtain ${{{\mathbf{\Psi}}}_{22}}=-\rho^*{{ {\mathbf{\widehat \Sigma}^{-1}} }}$. Then $\dot V =  - \beta{{\hat {\mathbf{x}} }^T}\widehat {\mathbf{\Gamma}} \left\{ \alpha \widehat {\mathbf{\Gamma}}  + {\rho ^*}{{ {\mathbf{\widehat\Sigma}^{-1}} }} \right\}\hat {\mathbf{x}}$, which is negative definite for any positive value of $\alpha$, $\beta$ and $\gamma$. Therefore, \eqref{eq:PID:1} achieves admissible consensus. Moreover, all nodes will converge to $x_\infty$ as defined in Prop. \ref{equil_point} with ${\sum\nolimits_{k = 1}^N {{\rho _k}}}=-N\rho^*$. To estimate the upper bound on the integral states, we consider \eqref{bound_z} with $\bar {{\boldsymbol\rho} }=\mathbb{0}_{1\times (N-1)}$ and using Proposition \ref{prop:Hmatr} we obtain expression \eqref{boundPI_heter}.
\end{IEEEproof}

\begin{Rem}
Note that:
\begin{itemize}
\item Network \eqref{eq:PID:1} still achieves admissible consensus if the nodes are unstable ($\rho^*\in\mathbb{R}^-$), since a suitable $\alpha$ can be chosen such that ${\alpha \widehat {\mathbf{\Gamma}}  + {\rho ^*}{\widehat{\mathbf{\Sigma}}^{-1}}}$ remains positive definite; however, the average trajectory ${{ {\hat x} }_1}(t)$ will be unstable and the system will exhibit unbounded admissible consensus; that is, ${\lim _{t \to \infty }}({x_i(t)-x_j(t)})=0$ but ${\lim _{t \to \infty }}{x_i}(t) \to \infty$ for all $i,j \in \mathcal{N}$.   

\item The advantage of the distributed PI action is to guarantee that all nodes converge towards the same asymptotic value at steady-state despite the presence of the disturbances. Therefore it should not be surprising that when these are present, consensus is indeed achieved but on a value that is dependent upon their magnitude.
\end{itemize}
\end{Rem}
%
%
%
We study next the distributed PI scheme obtained by setting $\gamma= 0$ in \eqref{eq:cont:2}. This case was already studied in the literature and a stability proof can also be found in \cite{Freeman2006} and \cite{Andreasson2012a} for identical and nonidentical disturbances, respectively. Nevertheless, for the sake of completeness, we report below an alternative proof based on our approach. 
\begin{coll}
\label{Coro:PI_homo}
Under the action of the distributed proportional-integral (PI) control obtained by setting  $\gamma= 0$ in \eqref{eq:cont:2}, a homogeneous network of $N$ linear systems \eqref{eq:sys:1}  achieves admissible consensus for any positive value of $\alpha$, and $\beta$. Moreover, all node states converge asymptotically to ${x_\infty}:=(1/N)\sum\nolimits_{k = 1}^N {{\delta _k}/} {\rho ^*}$ and ${z_\infty} \le \sqrt {N(N-1)} \left\| \mathbf{\Delta}  \right\|$ as $t$ approaches infinity.
\end{coll}
\begin{IEEEproof}
Setting $\gamma=0$ in \eqref{eq:PID:1}, yields $\tilde {\WideLaplacian}=\mathbf{I}_N$; therefore, $\tilde {\WideLaplacian}^{-1}=\mathbf{I}_N$ and their blocks in \eqref{inv_lapla} are $\widehat {l}_{11}=1$, $\widehat {\WideLaplacian}_{12}=\mathbb{0}$, $\widehat {\WideLaplacian}_{21}=\mathbb{0}$ and $\widehat {\WideLaplacian}_{22}=\mathbf{I}_{N-1}$. Then, ${\specnorm {\widehat{\mathbf{H}}}}=1$ in \eqref{eq:hat_H} and ${\mathbf{\Psi} _{22}} = N{\mathbf{R}_{22}}\left( {{\rho _1}{\mathbb{1}_{N - 1}}\mathbb{1}_{N - 1}^T + \widehat {\mathbf{P}}} \right)\mathbf{R}_{22}^T$ in \eqref{eq:PSI:22}. Since $\mathbf{\widehat P}=-\rho^* \mathbf{I}_{N-1}$ and, using properties (\ref{prop:U:2}) and (\ref{prop:U:4}), ${\mathbf{\Psi}_{22}} = -\rho^* \mathbf{I}_{N-1}$. Furthermore, setting $\gamma  = 0$ in \eqref{Gamma:1} we have ${\left.\widehat {\mathbf{\Gamma}}\right|_{\gamma  = 0}} = \widehat {\mathbf{\Lambda}}$
where
\begin{equation}
\label{eq:PI:De}
\widehat {\mathbf{\Lambda}} = \mbox{diag}\left\{ {{\lambda _2}, \cdots ,{\lambda _N}} \right\}
\end{equation}
It follows from the proof of Theorem \ref{TH:PID_integrators} that the Lyapunov function \eqref{eq:lyap} is still a positive, and globally unbounded, function. Thus, $\dot {V} =  -\beta {{\hat {\mathbf{x}} }^T}\widehat {\mathbf{\Lambda}} \left( {\alpha \widehat {\mathbf{\Lambda}}  + \rho^* \mathbf{I}_{N-1}} \right)\hat {\mathbf{x}}$, which is negative definite for any positive value of the control parameters $\alpha$, and $\beta$. Finally, we can conclude that the closed-loop network \eqref{eq:PID:1} achieves admissible consensus to ${x_\infty} :=   (1/N)\sum\nolimits_{k = 1}^N {{\delta _k}/} {\rho ^*}$ and substituting ${\specnorm {\widehat{\mathbf{H}}}}=1$ and $\bar {{\boldsymbol\rho} }=\mathbb{0}_{1\times (N-1)}$ in \eqref{bound_z} concludes the proof.
\end{IEEEproof}

By comparing the PID and PI strategies discussed above, we observe that the most notable difference is the presence of the factor $N/(\gamma\lambda_2+1)$ in the expression of the upper bound of the integral term when PID is used instead of PI. Such a factor can be varied by selecting the gain of the derivative action. This can be done by taking into account the size ($N$) and structure of the network encoded by $\lambda_2$, in order to avoid possible saturation of those integral terms and avoid the need for anti-windup strategies that can be difficult to implement across the network.
%
%
\begin{theor}(\textit{Convergence Rate})
\label{conve:rate:PIDcol}
The closed-loop network \eqref{eq:PID:1} with homogeneous node dynamics ($\rho_i=-\rho^*$), reaches admissible consensus with a convergence rate, say $\mu$, that can be estimated as
\begin{equation}
\label{conve:rate:PID} 
\mu=\left|\mathop {\max }\limits_{2 \le k \le N} \left(\mathbb{Re}\left(-\frac{b}{2}+\frac{\sqrt{b^2 -\frac{4\beta\lambda_k}{\gamma\lambda_k+1}}}{2} \right)\right)\right|
\end{equation} 
where $b={(\alpha\lambda_k+\rho*)}/({\gamma\lambda_k+1})$
\end{theor}
\begin{IEEEproof}
The proof is based on straightforward linear algebra \cite{Ren2007c}. Indeed, as all poles are identical, we have from \eqref{Psi:22:eq} that ${{{\mathbf{\Psi}}}_{22}}=-\rho^*{{ {\mathbf{\widehat \Sigma}^{-1}} }}$. Thus, one has from \eqref{eq:final:PID} that 
\[
\begin{array}{l}
\left[ {\begin{array}{*{20}{c}}
{\dot {\hat {\mathbf{x}}} (t)}\\
{\dot {\hat {\mathbf{z}}}(t)}
\end{array}} \right]=\underbrace{\left[ {\begin{array}{*{20}{c}}
{-\left( \rho^* \widehat{\mathbf{\Sigma}}^{-1} + \alpha \widehat {\mathbf{\Gamma}} \right)}&{\mathbf{I}_{N-1}}\\
{- \beta \widehat {\mathbf{\Gamma}}}&{\mathbb{0}}
\end{array}}\right]}_{\tilde{\mathbf{A}}}\left[ {\begin{array}{*{20}{c}}
{\hat {\mathbf{x}} (t)}\\
{\hat {\mathbf{z}} (t)}
\end{array}} \right]
\end{array}
\]

The rate of convergence can be estimated by computing the dominant eigenvalue(s) of the dynamic matrix $\tilde{\mathbf{A}}$ defined above. Specifically, say $\eta$ a generic eigenvalue of $\tilde{\mathbf{A}}$ and $\mathbf{v}=[\mathbf{v}_x\ \mathbf{v}_z]^T$ its corresponding eigenvector such that $\tilde{\mathbf{A}}\mathbf{v}=\eta\mathbf{v}$. Then, using the definition of $\tilde{\mathbf{A}}$ we obtain
\begin{subequations}
\begin{alignat}{3}
   \label{rate_conv_PIDa}
   -(\rho^*\widehat{\mathbf{\Sigma}}^{-1}+\alpha\widehat {\mathbf{\Gamma}})\mathbf{v}_x+\mathbf{v}_z &=& \eta\mathbf{v}_x  \\
    \label{rate_conv_PIDb}
   - \beta \widehat {\mathbf{\Gamma}}\mathbf{v}_x &=& \eta\mathbf{v}_z
\end{alignat}
\end{subequations}
Thus combining  \eqref{rate_conv_PIDb} and \eqref{rate_conv_PIDa}, we get
$$
-\left(\frac{\alpha\eta+\beta}{\eta}\right)\widehat {\mathbf{\Gamma}}\mathbf{v}_x=\left(\eta\mathbf{I}_{N-1}+\rho^*\widehat{\mathbf{\Sigma}}^{-1}\right)\mathbf{v}_x
$$
From their definitions, it is easy to see that all matrices on both sides are diagonal; hence, component-wise we obtain
$$
-\left(\alpha\eta+\beta\right)\frac{\lambda_k}{\gamma\lambda_k+1}=\eta\left(\eta+\frac{\rho^*}{\gamma\lambda_k+1}  \right), k\in\left\{2,\cdots,N\right\}
$$
Therefore, the eigenvalues of matrix $\tilde{\mathbf{A}}$ are the $2(N-1)$ solutions of the equations
$$
\eta_k^2+\eta_k\frac{(\alpha\lambda_k+\rho*)}{\gamma\lambda_k+1}+\frac{\beta\lambda_k}{\gamma\lambda_k+1} = 0, k\in\left\{2,\cdots,N\right\}
$$

Finally, letting $\mu:=\left|\eta_{max}\right|=\left|\mathop{max}\limits_{k}\left\{\mathbb{Re}\left(\eta_k^{\pm}\right)\right\}\right|$ we obtain \eqref{conve:rate:PID}. 
\end{IEEEproof}
\begin{Rem}
Note that
\begin{itemize}
\item The convergence rate depends on the network structure (via $\lambda_2$) as in the case of classical consensus problems, e.g. \cite{OlfatM2004}, but also on the node dynamics ($\rho*$), and the controller gains ($\alpha$, $\beta$ and $\gamma$). 

\item In general, increasing the value of $\gamma$ yields lower values of $b$ in \eqref{conve:rate:PID} and therefore the convergence rate may become slower. This indicates the presence of a trade-off between speed of convergence and bounds on the integral action that needs to be taken into account during the design stage. 
\end{itemize}
\end{Rem}

We conclude our investigation of the homogeneous case by studying the distributed proportional-derivative (PD) strategy obtained by setting $\beta=0$ in (\ref{eq:cont:2}).

\begin{theor}
\label{coll:PD}
Network \eqref{eq:PID:1} with homogeneous node dynamics, controlled by the distributed PD protocol obtained by selecting $\alpha>0$, $\beta=0$ and $\gamma>0$ in (\ref{eq:cont:2}) achieves $\varepsilon$-admissible consensus with 
\begin{equation}
\label{PD:bound}
\varepsilon=\frac{\gamma\lambda_N+1}{\gamma\lambda_2+1}\frac{N}{\alpha\lambda_N+\rho^*} \left\|\mathbf{\Delta}\right\| 
\end{equation}
\end{theor}
\begin{IEEEproof}
Equation \eqref{eq:cont:PID_Proj:a} without the integral action ($\beta=0$) and homogeneous nodes can be written as the two uncoupled equations
\begin{subequations}
\label{eq:DPD}
\begin{equation}
 {{\dot x }_1^\perp}(t) =   {{\psi _{11}}}{x _1^\perp}(t) + \mathbf{q}\mathbf{\Delta}
 \end{equation}
 \begin{equation}
 \label{PD:equ}
 \dot { {\bar{\mathbf{x}}}}^\perp (t) = \mathbf{D}{\mathbf{\bar x}}^\perp  + {\mathbf{\widehat R}}\mathbf{\Delta}
 \end{equation}
 \end{subequations}
where $\mathbf{D} = \left( {{ { {\mathbf{\Psi}} }_{22}} - \alpha  {\mathbf{\widehat \Gamma}} } \right)$. Note that, $\psi _{11}=-\rho^*$ and, using property (\ref{prop:L:5}) as done in the proof of Theorem \ref{TH:PID_integrators}, we have ${ {\mathbf{\Psi}} _{22}} =  - {\rho ^*}{{ {\mathbf{\widehat \Sigma}}^{-1} }}$ and $\mathbf{D}=-( \alpha \widehat {\mathbf{\Gamma}}  + {\rho ^*}{{ {\mathbf{\widehat \Sigma}^{ - 1}} }} )$ which is a negative definite, invertible matrix. Thus, ${\lim _{t \to \infty }}{ {\bar{\mathbf{x}}}}^\perp (t) = -\mathbf{D}^{-1}{\mathbf{\widehat R}}\mathbf{\Delta}$. 
To get an expression for the upper bound, we notice that, using \eqref{norm_R:a}, we can write
\[\begin{array}{l}
\specnorm{-\mathbf{D}^{-1}{\mathbf{\widehat R}}\mathbf{\Delta}}\leq  {\specnorm {\widehat{\mathbf{H}}}}\specnorm{( \alpha \widehat {\mathbf{\Gamma}}  + {\rho ^*}{{ {\mathbf{\widehat \Sigma}^{-1}} }} )^{-1} } \left\|\mathbf{\Delta}\right\|\end{array}\]
From \eqref{Gamma:1} and \eqref{bar_Sigma}, we know that $\widehat {\mathbf{\Gamma}}$ and $\widehat{\mathbf{\Sigma}}^{-1}$ are diagonal matrices and so is the matrix $\alpha\widehat {\mathbf{\Gamma}} + \rho^*\widehat{\mathbf{ \Sigma}}^{-1}$ with entries $(\alpha\lambda_k+\rho^*)/(\gamma\lambda_k+1),\forall k \in \{2,\cdots,N\}$. Therefore, the diagonal elements of $\left(\alpha\widehat {\mathbf{\Gamma}} + \rho^*\widehat{\mathbf{ \Sigma}}^{-1}\right)^{-1}$ are given by $(\gamma\lambda_k+1)/(\alpha\lambda_k+\rho^*)$. From Lemma \ref{lemm:quadratic_form}, we then have that $\specnorm{\left(\alpha\widehat {\mathbf{\Gamma}} + \rho^*\widehat{\mathbf{ \Sigma}}^{-1}\right)^{-1}} = (\gamma\lambda_N+1)/(\alpha\lambda_N+\rho^*)$ since $\lambda_N$ is the maximum eigenvalue of $\WideLaplacian$. Finally, using Prop. \ref{prop:Hmatr}, we obtain (\ref{PD:bound}).
\end{IEEEproof}
As expected, the bound $\varepsilon$ on the consensus error can be considerably reduced by increasing  the gain of the proportional action ($\alpha$) while it might be adversely affected by the gain of the derivative action. Indeed as for classical PID control, the presence of a distributed derivative action has little or no beneficial effect on the magnitude of the steady-state error. Also, it is clear that the network structure encoded by $\lambda_2$ and $\lambda_N$ has an effect on the overall error bound.
%
%
\subsection{Step 2b: Heterogeneous Node Dynamics}
Next we consider the case where at least one pole $\rho_i$ in \eqref{eq:PID:1} is different from the others, and the disturbances $\delta_i$ are generically nonidentical.
\begin{theor}
\label{TH:PID_integrators_heter}
The heterogeneous group of agents \eqref{eq:sys:1} controlled by the distributed PID strategy (\ref{eq:cont:2}), achieves admissible consensus for any $\beta>0$ and $\gamma\ge 0$ if the following conditions hold
\begin{subequations}
  \label{eq:cond:PID:heter}
	\begin{equation}
		\psi _{11} = (1/N)\sum\nolimits_{k = 1}^N {\rho _k}<0, \label{eq:cond:PID:heter:b}
	\end{equation}
	
\begin{equation}
\label{eq:cond:PID:heter:a}
\alpha{\frac{\lambda_2}{\gamma\lambda_2+1}} > \frac{1}{N}\left( {  \max_i\left\{ |\rho_i| \right\} + \frac{{\bar{\boldsymbol\rho}\bar{\boldsymbol\rho}^T}}{4\left| {{\psi _{11}}} \right| } \specnorm{\mathbf{H}_1}^2} \right)
\end{equation}
\end{subequations}
where $\mathbf{H}_1:=\mathbf{I}_{N-1} +  {\widehat{\mathbf{H}}}$.
Moreover, all node states converge to ${x_\infty}$ as defined in Prop. \ref{equil_point}, and the integral actions remain bounded by $z_\infty$ given in \eqref{bound_z}.
\end{theor}
\begin{IEEEproof}
Consider the following candidate Lyapunov function:
\begin{equation}
\label{Lyap:PI:hete}
V = \frac{1}{2}\left( {\hat x_1^2 + {{\hat {\mathbf{x}}}^T}\hat {\mathbf{x}} + \frac{1}{\beta }{{\hat {{\mathbf{z}}}}^T}\widehat {\mathbf{\Gamma}}^{-1} \hat {{\mathbf{z}}}} \right)
\end{equation}
where $\widehat {\mathbf{\Gamma}}^{-1}:=\mbox{diag}\{(\gamma\lambda_2+1)/\lambda_2,\cdots,(\gamma\lambda_N+1)/\lambda_N \}$. Differentiating $V$ along the trajectories of \eqref{eq:final:PID} yields 
\[
\begin{array}{l}
\dot {V} = \psi _{11}\hat x_1^2 + \hat {\mathbf{x}}^T{\mathbf{\Psi}} _{22}\hat {\mathbf{x}} - \alpha\hat {\mathbf{x}}^T{\widehat {\mathbf{\Gamma}} }\hat {\mathbf{x}} + 
\hat x_1\hat {\mathbf{x}}^T\left(\mathbf{\Psi} _{12}^T + {\mathbf{\Psi}} _{21}  \right)
\end{array}
\]
and, using \eqref{eq:PSI:12} and \eqref{eq:PSI:21}, we get
\begin{equation}
\label{diff:V:PID}
\begin{array}{l}
\dot {V} = \psi _{11}\hat x_1^2 + \hat {\mathbf{x}}^T{\mathbf{\Psi}} _{22}\hat {\mathbf{x}} - \alpha\hat {\mathbf{x}}^T{\widehat {\mathbf{\Gamma}} }\hat {\mathbf{x}} + 
\underbrace{\hat x_1\hat {\mathbf{x}}^T\mathbf{R}_{22}{\mathbf{H}_1}\bar{\boldsymbol\rho}^T}_{ g{\left(\hat x_1,\hat {\mathbf{x}}\right)}}
\end{array}
\end{equation}
Now, by setting $\mathbf{Q}^T=\mathbf{R}_{22}\mathbf{H}_1$, ${\boldsymbol\zeta}_1^T=\hat{\mathbf{x}}^T$, and ${\boldsymbol\zeta}_2=\hat x_1\bar{\boldsymbol\rho}^T$ in \eqref{pos:eq:1}, we can upper-bound $g{\left(\hat x_1,\hat {\mathbf{x}}\right)}$ as follows
\[
\begin{split}
g{\left(\hat x_1,\hat {\mathbf{x}}\right)} \le \frac{\sigma}{2}{\hat{\mathbf{x}}^T}\mathbf{Q}^T{\mathbf{Q}}\hat{\mathbf{x}} + \frac{\bar{\boldsymbol\rho}\bar{\boldsymbol\rho}^T}{2\sigma} \hat x_1^2
\end{split}
\]
This yields that from \eqref{diff:V:PID} we get
\begin{equation}
\label{diff:V:PID2}
\begin{array}{l}
\dot {V} \leq \left(\psi _{11} + \frac{\bar{\boldsymbol\rho}\bar{\boldsymbol\rho}^T}{2\sigma}  \right)\hat {x}_1^2 + \hat {\mathbf{x}}^T{\mathbf{\Psi}} _{22}\hat {\mathbf{x}} - \alpha\hat {\mathbf{x}}^T{\widehat {\mathbf{\Gamma}} }\hat {\mathbf{x}}   +\frac{\sigma}{2}{\hat{\mathbf{x}}^T}\mathbf{Q}^T{\mathbf{Q}}\hat{\mathbf{x}}
\end{array}
\end{equation}
From \eqref{eq:cont:Psieq}, we obtain $\specnorm{{\mathbf{\Psi}}+ {\mathbf{\Psi}}^T} \le 2\specnorm{\mathbf{U}}^2\specnorm{\mathbf{P}}\specnorm{\tilde{\WideLaplacian}^{-1}}$. Moreover, we have $\specnorm{\tilde{\WideLaplacian}^{-1}}=\left| \lambda_{max}\left( {\tilde{\WideLaplacian}^{-1}} \right)\right|=1$, $\specnorm{\mathbf{P}}=\max_i\{|\rho_i|\}$ and, from \eqref{prop:U:6}, $\specnorm{\mathbf{U}}^2=(1/N)$. Then, $\specnorm{{\mathbf{\Psi}}+ {\mathbf{\Psi}}^T} \le (2/N)\max_i\left\{ |\rho_i| \right\}$. 

Using Theorem 8.4.5 in \cite{Bernstein2009} we then find that $\lambda_{max}(\mathbf{\Psi}_{22}+\mathbf{\Psi}_{22}^T)\le \lambda_{max}(\mathbf{\Psi}+\mathbf{\Psi}^T) $ so that $\hat {\mathbf{x}}^T{\mathbf{\Psi}} _{22}\hat {\mathbf{x}}=(1/2)\hat {\mathbf{x}}^T\left({\mathbf{\Psi}} _{22}+{\mathbf{\Psi}} _{22}^T\right)\hat {\mathbf{x}}\le (1/N)\max_i\{|\rho_i|\}\hat {\mathbf{x}}^T\hat {\mathbf{x}}$.
Also, as $-\hat {\mathbf{x}}^T{\widehat {\mathbf{\Gamma}} }\hat {\mathbf{x}}\leq-{\lambda_2/(\gamma\lambda_2+1)}\hat {\mathbf{x}}^T\hat {\mathbf{x}}$, we obtain 
\begin{equation}
\label{diff:V:PIDineq}
\begin{array}{l}
\dot {V} \leq \left(\psi _{11} + \frac{\bar{\boldsymbol\rho}\bar{\boldsymbol\rho}^T}{2\sigma}  \right)\hat {x}_1^2 + \left({\frac{1}{N}\max_i\{|\rho_i|\} } \right.\\ \left. {-\alpha{\frac{\lambda_2}{\gamma\lambda_2+1}} + \frac{\sigma}{2N}\specnorm{\mathbf{H}_1}^2 } \right)\hat {\mathbf{x}}^T\hat {\mathbf{x}}
\end{array}
\end{equation}
%
Now, $\dot {V}$ is negative definite if the terms $\xi_1:=\psi _{11} + {\bar{\boldsymbol\rho}\bar{\boldsymbol\rho}^T}/({2\sigma}) $ and $\xi_2:=(1/N)\max_i\{|\rho_i|\}-\alpha{{\lambda_2}/({\gamma\lambda_2+1})} + ({\sigma}/{2N})\specnorm{\mathbf{H}_1}^2$ are both negative. 

From the assumptions, we have that $\psi _{11}<0$, therefore $\xi_1<0$ is ensured if we choose $\sigma>{\bar{\boldsymbol\rho}\bar{\boldsymbol\rho}^T}/(2\left| {{\psi _{11}}} \right|)$ [this is always possible as $\sigma$ is an arbitrary positive constant in \eqref{pos:eq:1}]. Then the condition $\xi_2<0$ can be fulfilled by selecting the control gains so as to satisfy \eqref{eq:cond:PID:heter:a}. 
Therefore, all agents in (\ref{eq:sys:1}) achieves admissible consensus to $x_\infty$ as defined in Prop. \ref{equil_point}, and the integral actions remain bounded by \eqref{bound_z} with $\hat{\specnorm{\mathbf{H}}}$ being bounded by \eqref{bound_H} which completes the proof.
\end{IEEEproof}
%
%
%
\section{Example} As a representative example, we consider the problem of achieving consensus in the network of $N$ linearized droop-controlled inverters which was  studied in \cite{Simpson-Porco2013}. The network equations are
\begin{subequations}
\label{eq:cont:PS}
\begin{alignat}{3}
\label{eq:cont:PSa}
{\dot \theta _i}(t) &= P_i^* - {P_i}(t)+v_i(t), \quad i \in \left\{ {1, \cdots ,N} \right\} \\
\label{eq:cont:PSb}
{P_i}(t) &= \sum\limits_{j = 1,j \ne i}^N {{E_i}{E_j}\left| {{Y_{ij}}} \right|\left( {{\theta _i} - {\theta _j}} \right)} 
\end{alignat}
\end{subequations}
where $\theta_i$ represents the phase of each inverter, $v_i(t)$ is the exogenous control signal, ${E_i}>0$ the nodal voltage, $Y_{ij}$ is the admittance between inverter $i$ and $j$. $P_i^*$ and $P_i(t)$ represent the normalized nominal active power injection  and the active electrical power exchanged with the other nodes, respectively. 

To achieve consensus we propose to use a combination of local and distributed control actions as also done in \cite{Andreasson2012a,Freeman2006} with the notable difference that we now consider the case of heterogeneous local state-feedback actions characterized by different gains which are deployed together with the  distributed PID strategy proposed in this paper.
In particular, we set
\begin{equation}
\label{PIDcontr:PWS}
{v_i}(t) = k_i\theta_i(t) + u_i(t)
\end{equation}
with $k_i$ being the local feedback gains and $u_i(t)$ the distributed PID action.
Then, letting $\omega_{ij}:={E_i}{E_j}\left| {{Y_{ij}}} \right|$ being the weights on the edges of the generator network and $\WideLaplacian$ the associated Laplacian matrix, the problem becomes that of proving convergence in the heterogeneous network given by
\begin{subequations}
\label{eq:exampl:1}
\begin{alignat}{3}
\label{eq:exampl:1a}
\dot {{\boldsymbol\theta}}(t) &= \tilde {\WideLaplacian}^{-1}\left( {\mathbf{K} - \tilde \alpha \WideLaplacian} \right){\boldsymbol\theta}(t) + \mathbf{z}(t) + \tilde {\WideLaplacian}^{-1}\mathbf{\Delta}\\
\label{eq:exampl:1b}
\dot {{\mathbf{z}}}(t) &=  - \beta \tilde {\WideLaplacian}^{-1}\WideLaplacian{\boldsymbol\theta}(t)
\end{alignat}
\end{subequations}
where ${\boldsymbol\theta}(t)=[{\theta _1}(t), \cdots ,{\theta _N}(t)]$,  $\tilde {\WideLaplacian}= \mathbf{I}_N + \gamma \WideLaplacian$, $\tilde \alpha=1+\alpha$, and $\mathbf{K}=\mbox{diag}\{k_1,\cdots,k_N\}$, $\mathbf{\Delta} := \mbox{diag}\{P_1^*,\cdots,P_N^*\}$. 
System \eqref{eq:exampl:1} has the same structure as \eqref{eq:PID:1}. Then, Theorem \ref{TH:PID_integrators_heter} can be used to tune the control gains and guarantee convergence in the case of a fixed network structure.
\begin{figure}[tbp]
\centering {
\subfigure[]
{\label{fig:PI_Fa}
{\includegraphics[scale=0.19]{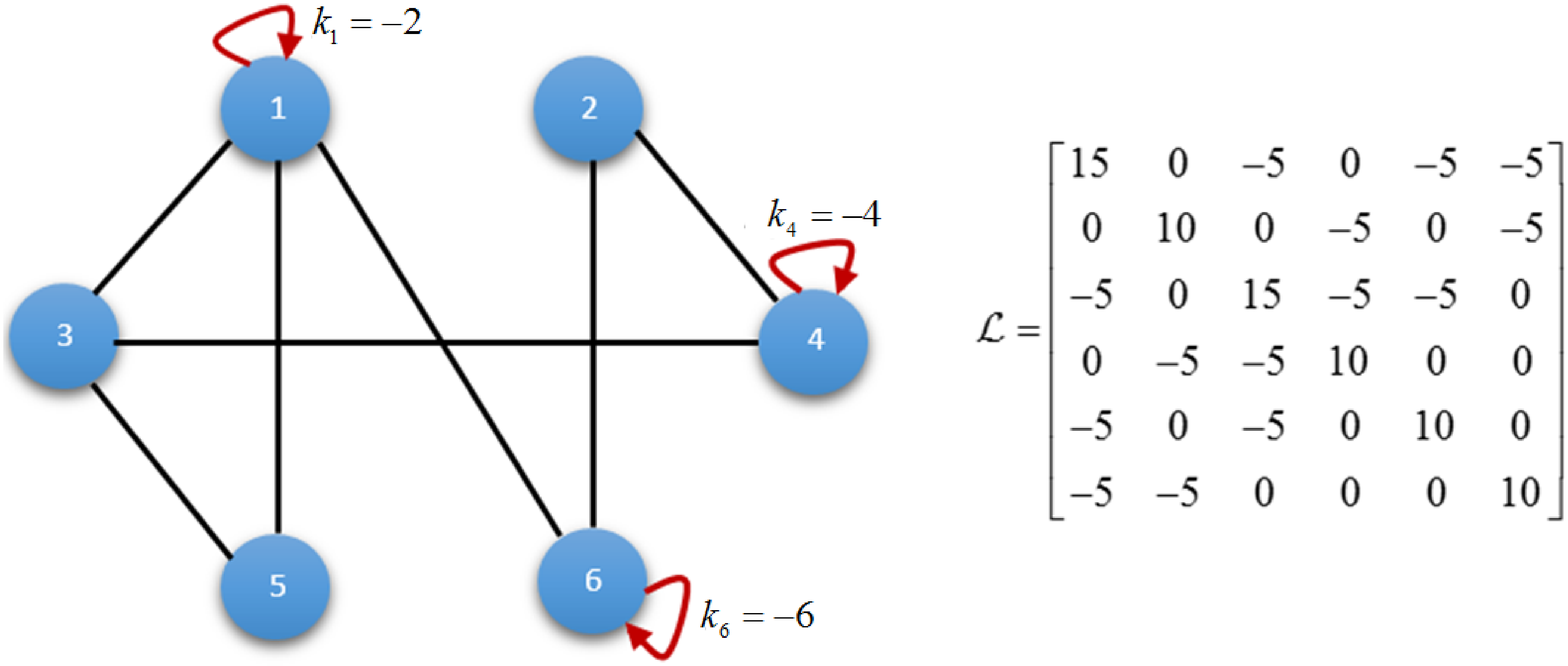}}}
\subfigure[]
{\label{fig:PI_Fb}
{\includegraphics[scale=0.17]{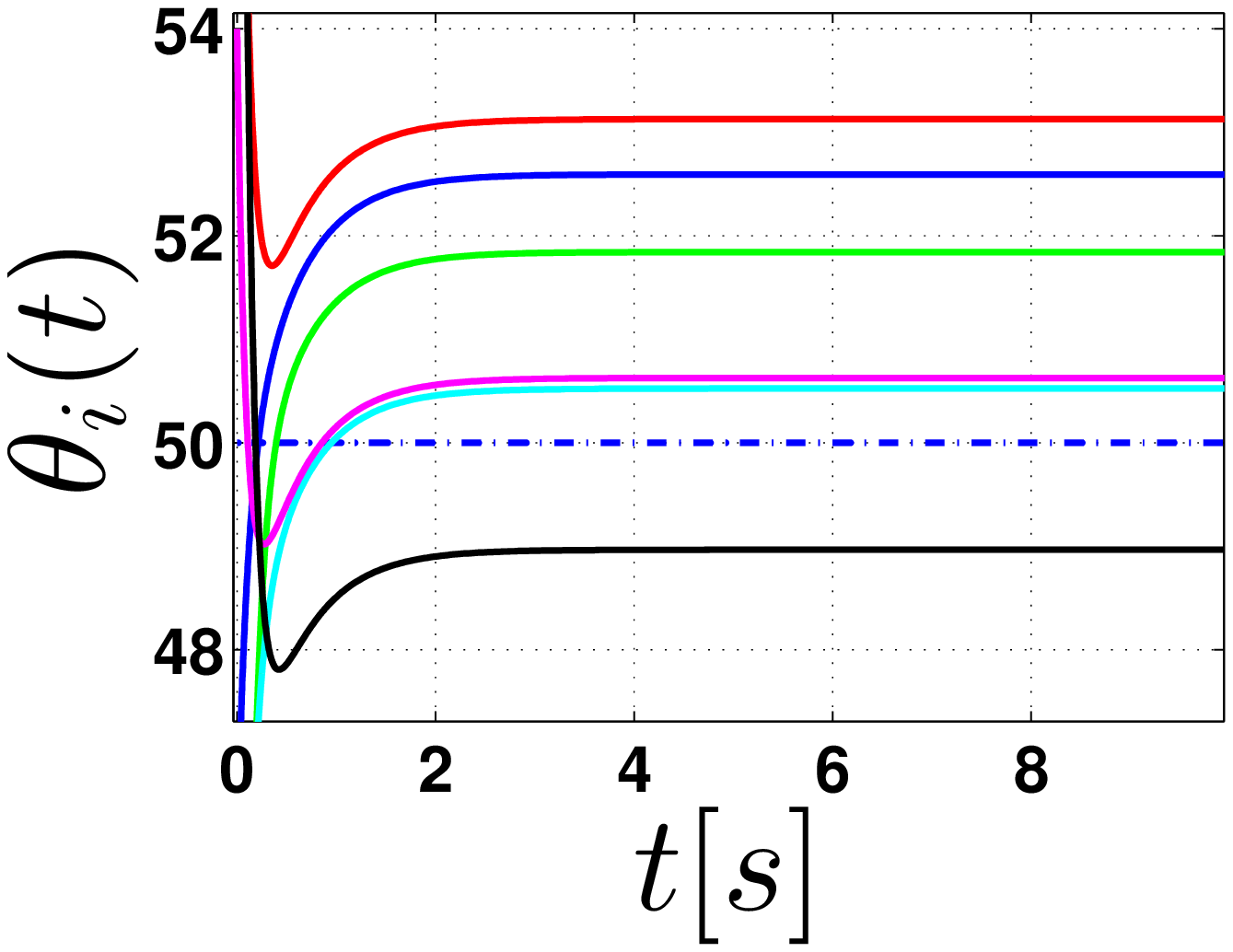}}}
\subfigure[]
{\label{fig:PI_Fd}
{\includegraphics[scale=0.17]{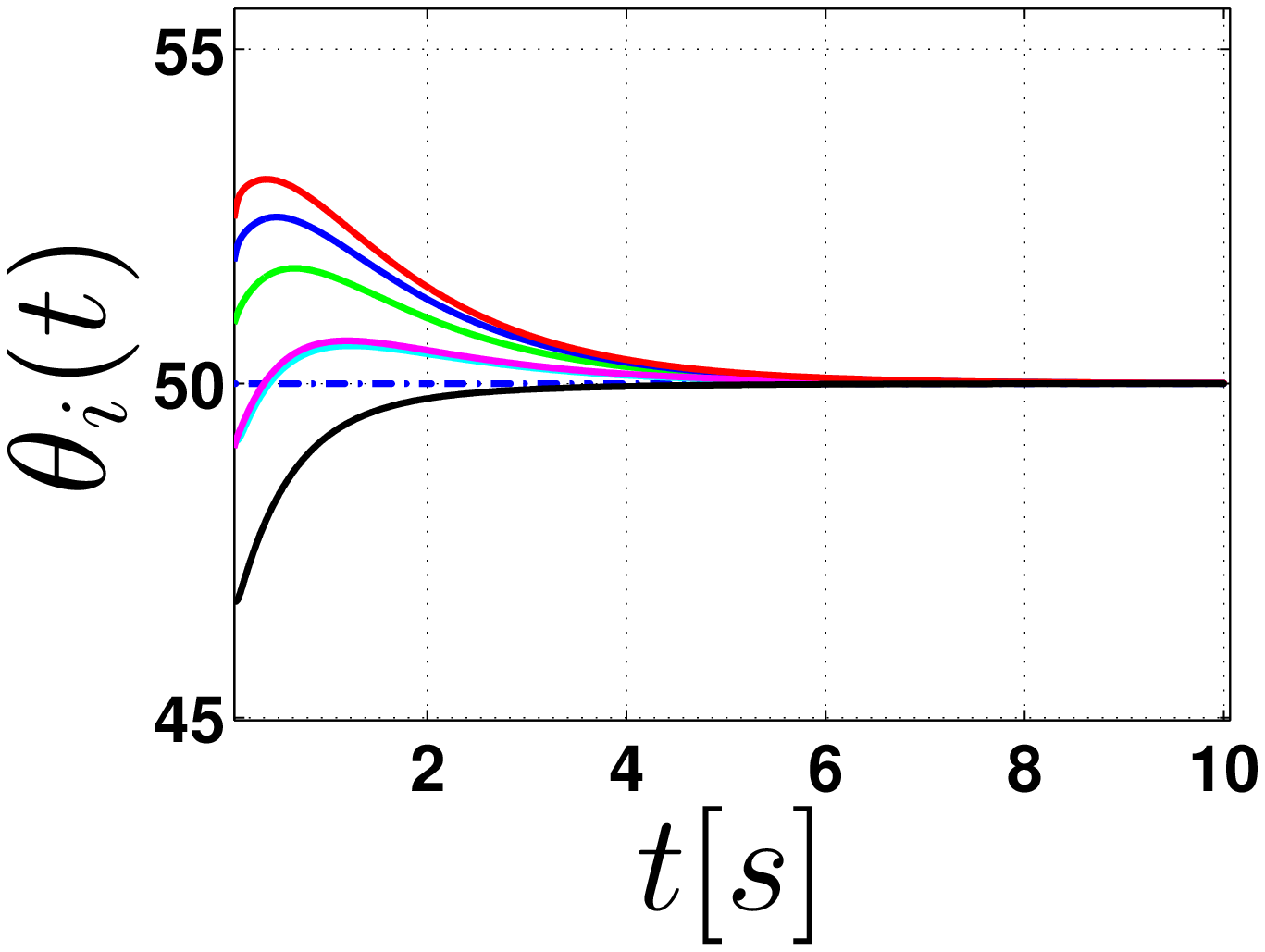}}}
\subfigure[]
{\label{fig:PI_Ff}
{\includegraphics[scale=0.17]{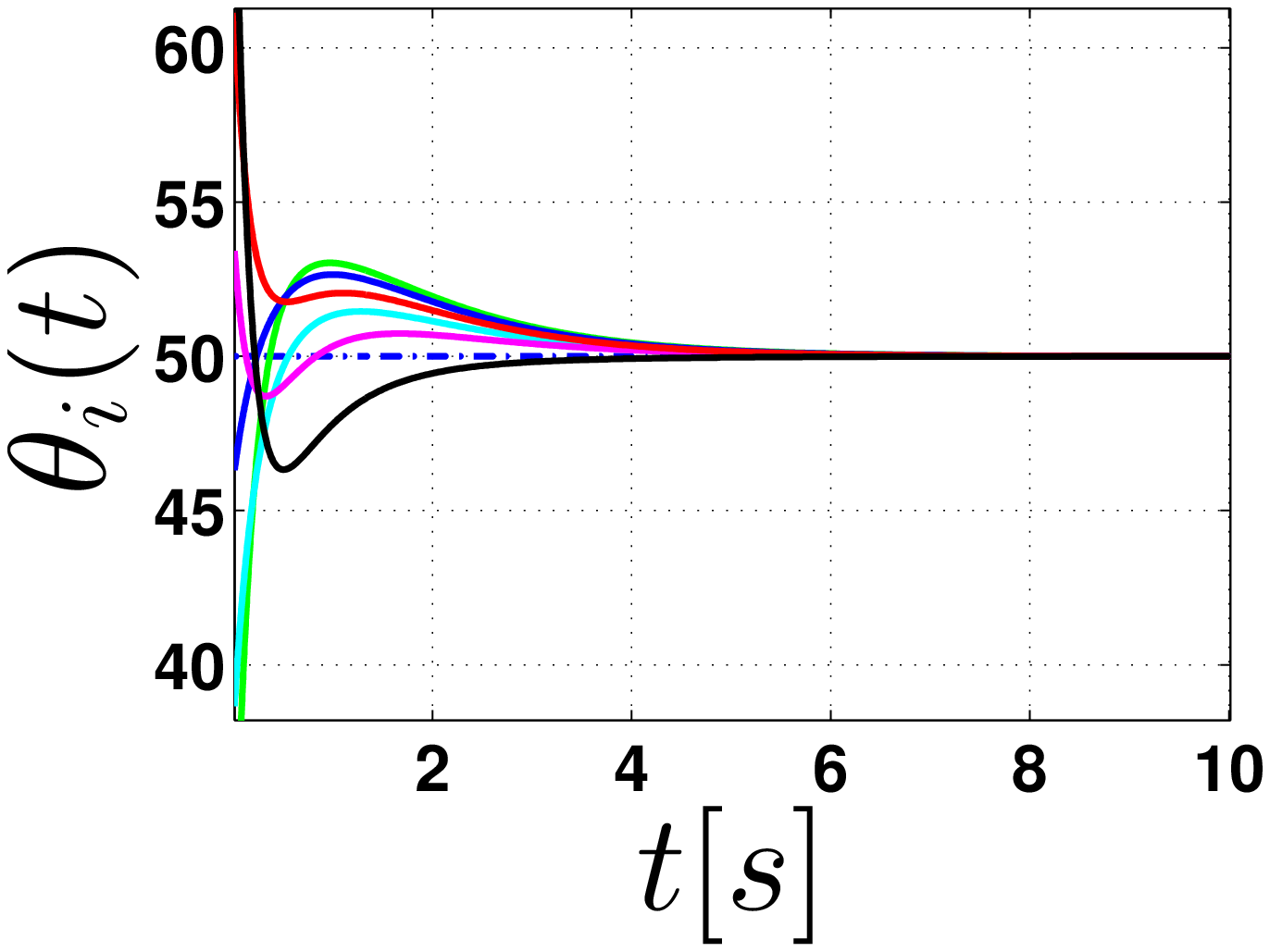}}}
\subfigure[]
{\label{fig:PI_Fc}
{\includegraphics[scale=0.17]{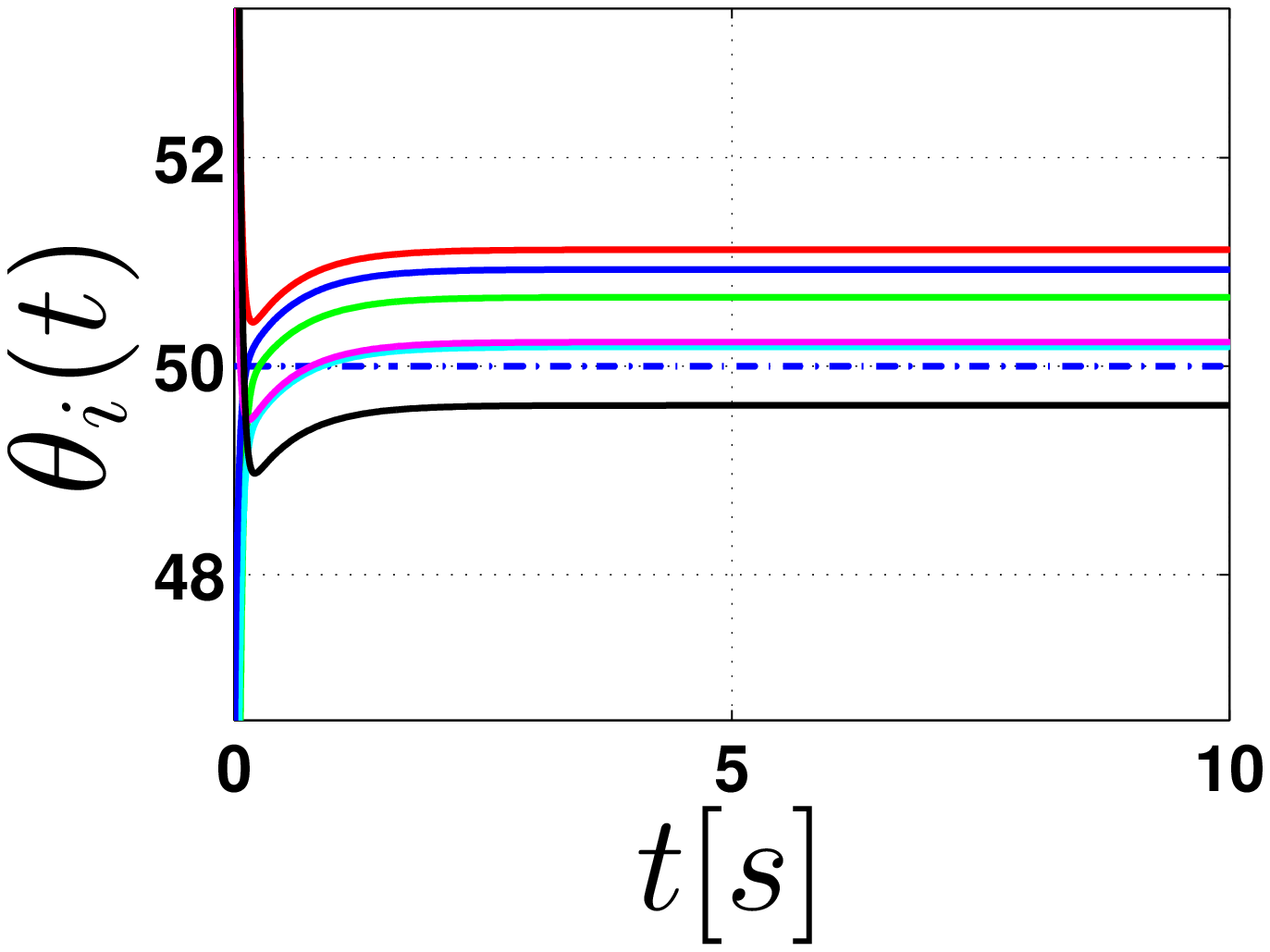}}}
\subfigure[]
{\label{fig:PI_Fe}
{\includegraphics[scale=0.17]{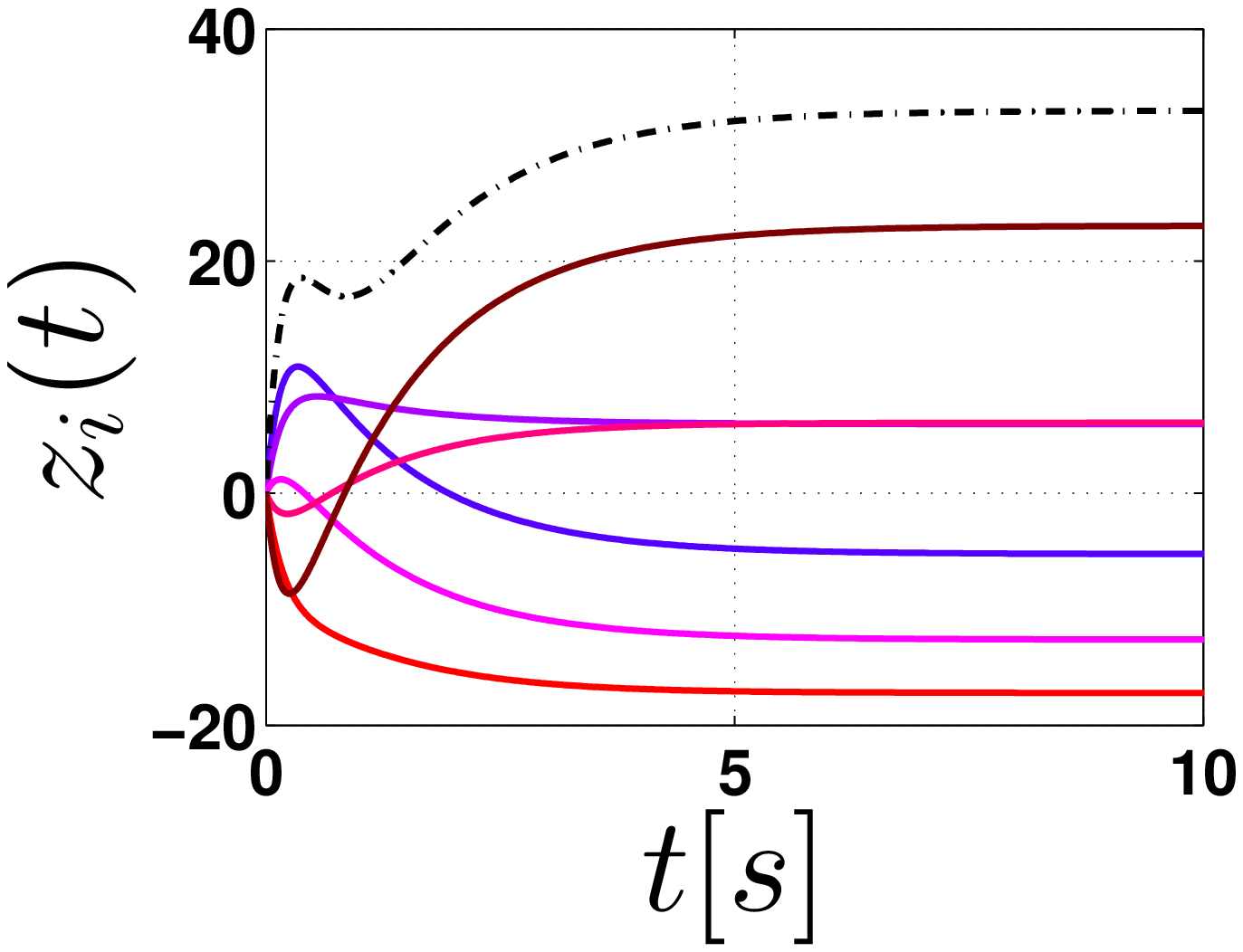}}}
\subfigure[]
{\label{fig:PI_Fg}
{\includegraphics[scale=0.17]{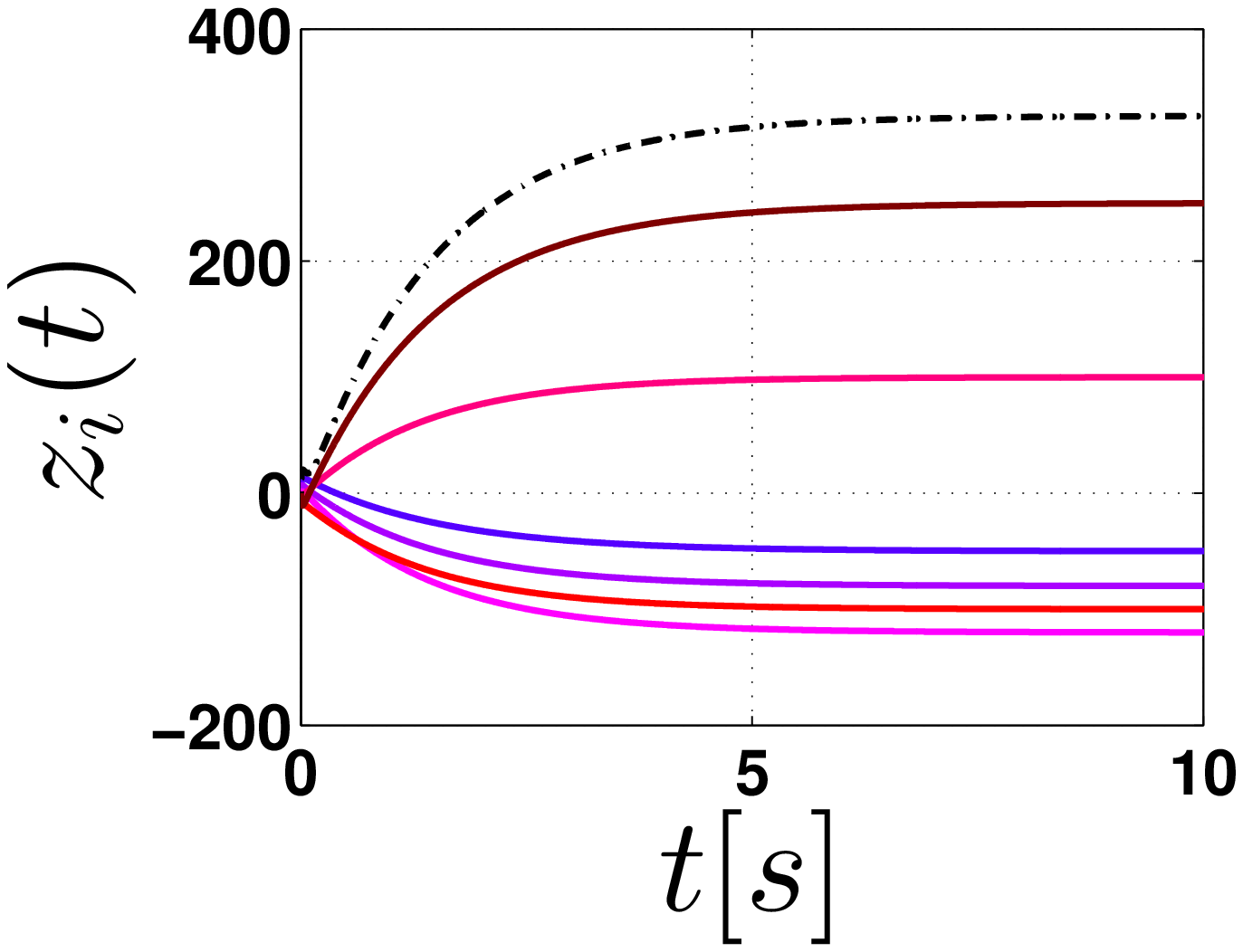}}}
}
\caption{(a) Schematic of a micro-grid of six inverters where all link weights are equal to 5. The red self-links represent local controllers acting on the node. Time response of the closed-loop network \eqref{eq:exampl:1} controlled by just a distributed proportional controller ($\beta=0$, $\gamma=0$) for (b) $\alpha=10$, (e) $\alpha=30$. (c),(f), and (d),(g) are the time evolution of the nodes and integral states for PID and PI respectively. The blue dash-dot line represent the theoretical convergence value $\theta_\infty$, and the black dash-dot line the time evolution of $\left\| {\mathbf{z}(t)} \right\|$. Self-loops in red represent the presence of a local feedback action.}
\label{fig:F_E} 
\end{figure}

As an illustration, consider the network shown in Fig. \ref{fig:PI_Fa} and assume the following nominal power injections in each node 
$\mathbf{\Delta}  = [150, 80, 120, 100, 100, 50]$.
According to what was reported in \cite{Kim2012}, purely proportional coupling leads to $\varepsilon$-admissible consensus. For example, in Figures \ref{fig:PI_Fb}, \ref{fig:PI_Fc}, we show the time response of the node dynamics in a heterogeneous network for two different values of $\alpha$ (with $\beta=\gamma=0$). We observe that in both cases a residual error is present that can only be reduced by increasing $\alpha$. To overcome these limitations, we consider now the same
network controlled via our strategy with $\mathbf{K} = \mbox{diag}\left\{-2,0,0,-4,0,-6\right\}$ and the gains of the distributed PID actions selected using Theorem \ref{TH:PID_integrators_heter}. 

Specifically, from the Laplacian matrix in Fig. \ref{fig:PI_Fa} we obtain $\lambda_2=5$. Also, we have $\psi_{11}=-2$ so that the first condition of Theorem \ref{TH:PID_integrators_heter} is fulfilled. Moreover, $\bar {\boldsymbol\rho} {\bar {\boldsymbol\rho} }^T=32$ and selecting $\gamma=1$ we obtain $\alpha > 5.92$. Without loss of generality, we choose $\alpha=6$, and $\beta=5$.
The resulting evolution of the node states and integral actions is shown in Fig. \ref{fig:PI_Fd}, \ref{fig:PI_Fe}, where admissible consensus is reached as expected to the predicted value $\theta_\infty:= - \sum\nolimits_{i = 1}^N {{\delta _i}} /\sum\nolimits_{i = 1}^N {{k _i}} = 50$. For the sake of comparison, the performance of a purely distributed proportional and integral action obtained by switching off the derivative actions in the previous example ($\gamma=0$) is depicted in Fig. \ref{fig:PI_Ff}, \ref{fig:PI_Fg}. We observe that the distributed PID strategy presented in this paper is indeed able to  guarantee better performance since when PID is used the bounds on the integral terms are smaller than those when controlled by a distributed PI strategy.
\section{Conclusions and Future Work}
We have investigated the use of  a distributed PID protocol to achieve consensus in homogeneous and heterogeneous multi-agent networks. Convergence of the strategy in both cases was obtained by using appropriate state transformations, linear algebra and Lyapunov functions. Explicit expressions for the consensus values were obtained together with analytical estimates of the upper bound for all integral actions. Also, some estimates of the rate of convergence were obtained as functions of the gains of the distributed control actions and the network structure. It was found that the network architecture, the nodal dynamics and the control gains all contribute to determine the stability and performance of the closed-loop network. 
\bibliographystyle{IEEEtran}
\bibliography{Distributed_PID_control.bbl}
\appendices
\section{Derivation of Expression \eqref{prop:L:5}}
\label{Appendix_I}
From Lemma \ref{lem:LDL_eig} we have that $\mathbf{U}^{ - 1}\tilde {\WideLaplacian}^{-1}\mathbf{U} = {\mathbf{\Sigma} ^{ - 1}} $. Then Using the block representation we have

\[
 \left[ {\begin{array}{*{20}{c}}
   r_{11} & {{\mathbf{R}_{12}}}  \\
   {\mathbf{R}_{21}} & {{\mathbf{R}_{22}}}  \\
\end{array}} \right]  \left[ {\begin{array}{*{20}{l}}
{\widehat {l}_{11}}&{\widehat {\WideLaplacian}_{12}}\\
{\widehat {\WideLaplacian}_{21}}&{\widehat {\WideLaplacian}_{22}}
\end{array}} \right] \left[ {\begin{array}{*{20}{c}}
   1 & {N{\mathbf{R}_{21}^T}}  \\
   {{\mathbb{1}_{N - 1}}} & {N{\mathbf{R}_{22}^T}}  \\
\end{array}} \right] = {\mathbf{\Sigma} ^{ - 1}}
\] 

Letting $\mathbf{M}:=\mathbf{U}^{ - 1}\tilde {\WideLaplacian}\mathbf{U}$, some straightforward algebra yields
\[\begin{array}{l}
\left[ {\begin{array}{*{20}{c}}
{\mathbf{M}_{11}}  &   {\mathbf{M}_{12}}\\
{\mathbf{M}_{21}} & {\mathbf{M}_{22}}
\end{array}} \right]
 = \left[ {\begin{array}{*{20}{c}}
1&\mathbb{0}_{1 \times (N - 1)}\\
\mathbb{0}_{(N - 1) \times 1}&{{{\mathbf{\widehat \Sigma} }^{ - 1}}}
\end{array}} \right]
\end{array}\]
where
\[
\begin{split}
\mathbf{M}_{11}& = {r_{11}}\left( {\widehat l_{11} + \widehat {\WideLaplacian}_{12}{\mathbb{1}_{N - 1}}} \right) + {\mathbf{R}_{12}}\left( {\widehat {\WideLaplacian}_{21} + \widehat {\WideLaplacian}_{22}{\mathbb{1}_{N - 1}}} \right)\\
\mathbf{M}_{12}& =  N{r_{11}}\left( {\widehat l_{11}\mathbf{R}_{21}^T + \widehat {\WideLaplacian}_{12}\mathbf{R}_{22}^T} \right)\\ 
& \quad +N{\mathbf{R}_{12}}\left( {\widehat {\WideLaplacian}_{21}{\mathbf{R}_{21}^T} + \widehat {\WideLaplacian}_{22}\mathbf{R}_{22}^T} \right)\\
\mathbf{M}_{21} &= {\mathbf{R}_{21}}\left( {\widehat l_{11} + \widehat {\WideLaplacian}_{12}{\mathbb{1}_{N - 1}}} \right) + {\mathbf{R}_{22}}\left( {\widehat {\WideLaplacian}_{21} + \widehat {\WideLaplacian}_{22}{\mathbb{1}_{N - 1}}} \right)\\
\mathbf{M}_{22} &= N{\mathbf{R}_{21}}\left( {\widehat l_{11}\mathbf{R}_{21}^T + \widehat {\WideLaplacian}_{12}\mathbf{R}_{22}^T} \right)\\
& \quad +N{\mathbf{R}_{22}}\left( {\widehat {\WideLaplacian}_{21}\mathbf{R}_{21}^T + \widehat {\WideLaplacian}_{22}\mathbf{R}_{22}^T} \right)
\end{split}
\]
Equating the blocks we have that $\mathbf{M}_{22}=\mathbf{\widehat \Sigma}^{ - 1}$, and some algebraic manipulations yield ${\mathbf{R}_{21}}\widehat l_{11}\mathbf{R}_{21}^T + {\mathbf{R}_{21}}\widehat {\WideLaplacian}_{12}\mathbf{R}_{22}^T + {\mathbf{R}_{22}}\widehat {\WideLaplacian}_{21}\mathbf{R}_{21}^T + {\mathbf{R}_{22}}\widehat {\WideLaplacian}_{22}\mathbf{R}_{22}^T = \frac{1}{N}{{ \mathbf{\mathbf{\widehat \Sigma}} }^{ - 1}}$. Now, adding and subtracting $\widehat l_{11}{\mathbf{R}_{21}}\mathbb{1}_{N - 1}^T\mathbf{R}_{22}^T$ one gets $\widehat l_{11}{\mathbf{R}_{21}}\left( {\mathbf{R}_{21}^T + \mathbb{1}_{N - 1}^T\mathbf{R}_{22}^T} \right) - \widehat l_{11}{\mathbf{R}_{21}}\mathbb{1}_{N - 1}^T\mathbf{R}_{22}^T + {\mathbf{R}_{21}}\widehat {\WideLaplacian}_{12}\mathbf{R}_{22}^T + {\mathbf{R}_{22}}\widehat {\WideLaplacian}_{21}\mathbf{R}_{21}^T + {\mathbf{R}_{22}}\widehat {\WideLaplacian}_{22}\mathbf{R}_{22}^T = \frac{1}{N}{{\widehat {\mathbf{\Sigma}} }^{ - 1}}$. From property (\ref{prop:U:3}) one has that $\mathbf{R}_{21}^T + {\mathbb{1}_{N - 1}}\mathbf{R}_{22}^T = \mathbb{0}$. Also, using (\ref{prop:U:1}), we have ${\mathbf{R}_{21}} =  - {\mathbf{R}_{22}}{\mathbb{1}_{N - 1}}$ so that the equation above can be recast as $ {\mathbf{R}_{22}}\widehat {\WideLaplacian}_{22}\mathbf{R}_{22}^T - {\mathbf{R}_{22}}{\mathbb{1}_{N - 1}}\widehat {\WideLaplacian}_{12}\mathbf{R}_{22}^T - {\mathbf{R}_{22}}\widehat {\WideLaplacian}_{21}\mathbb{1}_{N - 1}^T\mathbf{R}_{22}^T + \widehat l_{11}{\mathbf{R}_{22}}\mathbb{1}_{N - 1}^T{\mathbb{1}_{N - 1}}\mathbf{R}_{22}^T = \frac{1}{N}{{\mathbf{\widehat \Sigma} }^{ - 1}}$. Finally regrouping terms we obtain \eqref{prop:L:5}.  
\section{Computation of $\mathbf{ \Psi}$ matrix}
\label{Appendix_II}

We know that $\mathbf{ \Psi}  = {\mathbf{U}^{ - 1}}{{\tilde {\WideLaplacian}}^{ - 1}} \mathbf{P} \mathbf{U}$, and 
\[
\left[ {\begin{array}{*{20}{c}}
\psi_{11} &{{\mathbf{\Psi} _{12}}}\\
{{{ \mathbf{\Psi} }_{12}}}&{{{\mathbf{ \Psi} }_{22}}}
\end{array}} \right] = {\mathbf{U}^{ - 1}}{{\tilde {\WideLaplacian}}^{ - 1}} \left[ {\begin{array}{*{20}{c}}
1&\mathbb{0}_{1 \times (N - 1)}\\
\mathbb{0}_{(N - 1) \times 1}&{{{\mathbf{\widehat {\mathbf{P}}} }}}
\end{array}} \right] \mathbf{U} 
\]

Next, we simplify the expression of each block in the matrix. Specifically, from \eqref{eq:blockdef} we have that $\mathbf{R}_{12}=r_{11}\mathbb{1}_{N-1}^T$ and the first block can be expressed as
\[
\begin{split}
{\psi _{11}} &= {r_{11}}\left( {\widehat l_{11}{\rho _1} + \widehat {\WideLaplacian}_{12}\mathbf{\widehat{\mathbf{P}}}{\mathbb{1}_{N - 1}}} \right) \\ & \quad + {\mathbf{R}_{12}}\left( {\widehat {\WideLaplacian}_{21}{\rho _1} + \widehat {\WideLaplacian}_{22}\mathbf{\widehat {\mathbf{P}}}{\mathbb{1}_{N - 1}}} \right)\\
& = \,{r_{11}}\widehat l_{11}{\rho _1} + {r_{11}}\mathbb{1}_{N - 1}^T\widehat {\WideLaplacian}_{21}{\rho _1} + {r_{11}}\widehat {\WideLaplacian}_{12}\mathbf{\widehat {\mathbf{P}}}{\mathbb{1}_{N - 1}} \\
& \quad + {\mathbf{R}_{12}}\widehat {\WideLaplacian}_{22}\mathbf{\widehat {\mathbf{P}}}{\mathbb{1}_{N - 1}}
\end{split}
\]
From (\ref{prop:L:1}) one has $\mathbb{1}_{N - 1}^T\widehat {\WideLaplacian}_{21} = 1 - \widehat l_{11}$, where $\widehat {\WideLaplacian}_{12}=\widehat {\WideLaplacian}_{21}^T$. Thus, using (\ref{prop:L:2}) yields
\[
\begin{split}
{\psi _{11}} &= {r_{11}}{\rho _1} + {r_{11}}\widehat {\WideLaplacian}_{12}\mathbf{\widehat {\mathbf{P}}}{\mathbb{1}_{N - 1}} + {r_{11}}\mathbb{1}_{N - 1}^T\widehat {\WideLaplacian}_{22}\mathbf{\widehat {\mathbf{P}}}{\mathbb{1}_{N - 1}}\\
&= {r_{11}}{\rho _1} + {r_{11}}\widehat {\WideLaplacian}_{12}\mathbf{\widehat{\mathbf{P}}}{\mathbb{1}_{N - 1}}+{r_{11}}(\mathbb{1}_{N - 1}^T - \widehat {\WideLaplacian}_{12})\widehat {\mathbf{P}}{\mathbb{1}_{N - 1}}\\
&= {r_{11}}{\rho _1} + {r_{11}}\mathbb{1}_{N - 1}^T\mathbf{\widehat {\mathbf{P}}}{\mathbb{1}_{N - 1}}\\
&= (1/N)\sum\nolimits_{k = 1}^N {{\rho _k}} 
\end{split}
\]
We move next to the second block  given by  $\mathbf{\Psi}_{12} = N{r_{11}}( {\widehat l_{11}{\rho _1}\mathbf{R}_{21}^T + \widehat {\WideLaplacian}_{12}\mathbf{\widehat {\mathbf{P}}}\mathbf{R}_{22}^T}) + N{\mathbf{R}_{12}}( {{\rho _1}\widehat {\WideLaplacian}_{21}\mathbf{R}_{21}^T + \widehat {\WideLaplacian}_{22}\mathbf{\widehat {\mathbf{P}}}\mathbf{R}_{22}^T} ) $. From \eqref{eq:blockdef} we have that $\mathbf{R}_{12}=r_{11}\mathbb{1}_{N-1}^T$. Some algebraic manipulation yields $(1/N){\mathbf{\Psi} _{12}} = {r_{11}}{\rho _1}\widehat l_{11}\mathbf{R}_{21}^T + {\rho _1}{r_{11}}\mathbb{1}_{N - 1}^T\widehat {\WideLaplacian}_{21}\mathbf{R}_{21}^T + {r_{11}}\widehat {\WideLaplacian}_{12}\mathbf{\widehat {\mathbf{P}}}\mathbf{R}_{22}^T + {\mathbf{R}_{12}}\widehat {\WideLaplacian}_{22}\mathbf{\widehat {\mathbf{P}}}\mathbf{R}_{22}^T$.
Using (\ref{prop:L:1}) the right-hand side of this expression can be rewritten so as to get
\[
\begin{split}
\frac{1}{N}{\mathbf{\Psi} _{12}} &= {r_{11}}{\rho _1}\widehat l_{11}\mathbf{R}_{21}^T + {\rho _1}{r_{11}}(1 - \widehat l_{11})\mathbf{R}_{21}^T \\
& \quad+ {r_{11}}\widehat {\WideLaplacian}_{12}\mathbf{\widehat {\mathbf{P}}}\mathbf{R}_{22}^T + {\mathbf{R}_{12}}\widehat {\WideLaplacian}_{22}\mathbf{\widehat {\mathbf{P}}}\mathbf{R}_{22}^T\\
 \frac{1}{N}{\mathbf{\Psi} _{12}} &= {r_{11}}{\rho _1}\mathbf{R}_{21}^T + {r_{11}}\widehat {\WideLaplacian}_{12}\mathbf{\widehat {\mathbf{P}}}\mathbf{R}_{22}^T + {r_{11}}\mathbb{1}_{N - 1}^T\widehat {\WideLaplacian}_{22}\mathbf{\widehat {\mathbf{P}}}\mathbf{R}_{22}^T\\
 \frac{1}{N}{\mathbf{\Psi} _{12}} &= {r_{11}}{\rho _1}\mathbf{R}_{21}^T \\ 
 & \quad + {r_{11}}\widehat {\WideLaplacian}_{12}\mathbf{\widehat {\mathbf{P}}}\mathbf{R}_{22}^T + {r_{11}}(\mathbb{1}_{N - 1}^T - \widehat {\WideLaplacian}_{12})\mathbf{\widehat {\mathbf{P}}}\mathbf{R}_{22}^T\\
 \frac{1}{N}{\mathbf{\Psi} _{12}} &= {r_{11}}{\rho _1}\mathbf{R}_{21}^T + {\mathbf{R}_{12}}\mathbf{\widehat {\mathbf{P}}}\mathbf{R}_{22}^T
\end{split}
\]
Finally, adding and subtracting $\rho _1\mathbf{R}_{12}\mathbf{R}_{22}^T$ and applying property (\ref{prop:U:3}) yields
\[
\begin{split}
(1/N){\mathbf{\Psi} _{12}} &=  - {\rho _1}{\mathbf{R}_{12}}\mathbf{R}_{22}^T + {\mathbf{R}_{12}}\mathbf{\widehat P}\mathbf{R}_{22}^T\\
(1/N){\mathbf{\Psi}_{12}} &= {\mathbf{R}_{12}}\left( {\mathbf{\widehat P} - {\rho _1}{\mathbf{I}_{N - 1}}} \right)\mathbf{R}_{22}^T\\
(1/N){\mathbf{\Psi} _{12}} &= {r_{11}}\mathbb{1}_{N - 1}^T\left( {\mathbf{\widehat P} - {\rho _1}{\mathbf{I}_{N - 1}}} \right)\mathbf{R}_{22}^T\\
{\mathbf{\Psi}_{12}} &= [{\rho _2} - {\rho _1}, \cdots ,{\rho _N} - {\rho _1}]\mathbf{R}_{22}^T  = \bar {\boldsymbol\rho} \mathbf{R}_{22}^T
\end{split}
\]
The third block of the matrix has an expression given by ${{\mathbf{\Psi} }_{21}} = {\mathbf{R}_{21}}\widehat l_{11}{\rho _1} + {\mathbf{R}_{22}}\widehat {\WideLaplacian}_{21}{\rho _1} + {\mathbf{R}_{21}}\widehat {\WideLaplacian}_{12}\mathbf{\widehat {\mathbf{P}}}{\mathbb{1}_{N - 1}} + {\mathbf{R}_{22}}\widehat {\WideLaplacian}_{22}\mathbf{\widehat {\mathbf{P}}}{\mathbb{1}_{N - 1}}$. Using (\ref{prop:U:1}) we get
\[
\begin{split}
{{\mathbf{ \Psi} }_{21}} &=  - {\mathbf{R}_{22}}{\mathbb{1}_{N - 1}}\widehat l_{11}{\rho _1} + {\mathbf{R}_{22}}\widehat {\WideLaplacian}_{21}{\rho _1} \\
& \quad + {\mathbf{R}_{22}}\widehat {\WideLaplacian}_{22}\mathbf{\widehat{\mathbf{P}}}{\mathbb{1}_{N - 1}} - {\mathbf{R}_{22}}{\mathbb{1}_{N - 1}}\widehat {\WideLaplacian}_{12}\mathbf{\widehat {\mathbf{P}}}{\mathbb{1}_{N - 1}}\\
{{\mathbf{ \Psi} }_{21}} &= {\mathbf{R}_{22}}( ( {\widehat {\WideLaplacian}_{22} - {\mathbb{1}_{N - 1}}\widehat {\WideLaplacian}_{12}} )\mathbf{\widehat{\mathbf{P}}}{\mathbb{1}_{N - 1}} \\
 & \quad+ {\rho _1}( {\widehat {\WideLaplacian}_{21} - {\mathbb{1}_{N - 1}}\widehat l_{11}} ) )
\end{split}
\]
Then, using \eqref{prop:L:4}, one finally has ${{\mathbf{\Psi}}_{21}} = {\mathbf{R}_{22}}( {\widehat {\WideLaplacian}_{22} - {\mathbb{1}_{N - 1}}\widehat {\WideLaplacian}_{12}} )( {\mathbf{\widehat{\mathbf{P}}} - {\rho _1}{\mathbf{I}_{N - 1}}}){\mathbb{1}_{N - 1}}$ and \eqref{eq:PSI:21} is then obtained.
Finally, the last block of $\mathbf{\Psi}$ matrix is expressed as
\[
\begin{split}
(1/N){{\mathbf{\Psi} }_{22}} &= {\mathbf{R}_{21}}\widehat l_{11}{\rho _1}\mathbf{R}_{21}^T + {\mathbf{R}_{22}}{\rho _1}\widehat {\WideLaplacian}_{21}\mathbf{R}_{21}^T \\ 
& \quad + {\mathbf{R}_{22}}\widehat {\WideLaplacian}_{22}\mathbf{\widehat{\mathbf{P}}}\mathbf{R}_{22}^T + {\mathbf{R}_{21}}\widehat {\WideLaplacian}_{12}\mathbf{\widehat {\mathbf{P}}}\mathbf{R}_{22}^T
\end{split}
\]
Using (\ref{prop:U:1}) and (\ref{prop:L:1}) one gets
\[
\begin{split}
(1/N){{\mathbf{ \Psi} }_{22}} &= {\rho _1}{\mathbf{R}_{22}}( {\widehat {\WideLaplacian}_{21} - \widehat l_{11}{\mathbb{1}_{N - 1}}} )\mathbf{R}_{21}^T\\ 
& \quad + {\mathbf{R}_{22}}( {\widehat {\WideLaplacian}_{22} - {\mathbb{1}_{N - 1}}\widehat {\WideLaplacian}_{21}} )\mathbf{\widehat{\mathbf{P}}}\mathbf{R}_{22}^T\\
 &= {\mathbf{R}_{22}}( {\widehat {\WideLaplacian}_{22} - {\mathbb{1}_{N - 1}}\widehat {\WideLaplacian}_{12}} )\mathbf{\widehat {\mathbf{P}}}\mathbf{R}_{22}^T \\
& \quad - {\rho _1}{\mathbf{R}_{22}}( {\widehat {\WideLaplacian}_{21} - \widehat l_{11}{\mathbb{1}_{N - 1}}} )\mathbb{1}_{N-1}^T\mathbf{R}_{22}^T
\end{split}
\]
then, applying (\ref{prop:L:4}) one has
\[
\begin{split}
(1/N){{\mathbf{ \Psi} }_{22}} &= {\mathbf{R}_{22}}( {\widehat {\WideLaplacian}_{22} - {\mathbb{1}_{N - 1}}\widehat {\WideLaplacian}_{12}} )\mathbf{\widehat {\mathbf{P}}}\mathbf{R}_{22}^T\\ 
& \quad + {\rho _1}{\mathbf{R}_{22}}( {\widehat {\WideLaplacian}_{22} - {\mathbb{1}_{N - 1}}\widehat {\WideLaplacian}_{12}} ){\mathbb{1}_{N - 1}}\mathbb{1}_{N - 1}^T\mathbf{R}_{22}^T
\end{split}
\]
and we obtain \eqref{eq:PSI:22}
%
%
\begin{IEEEbiography}[{\includegraphics[width=1in,height=1.25in,clip,keepaspectratio]{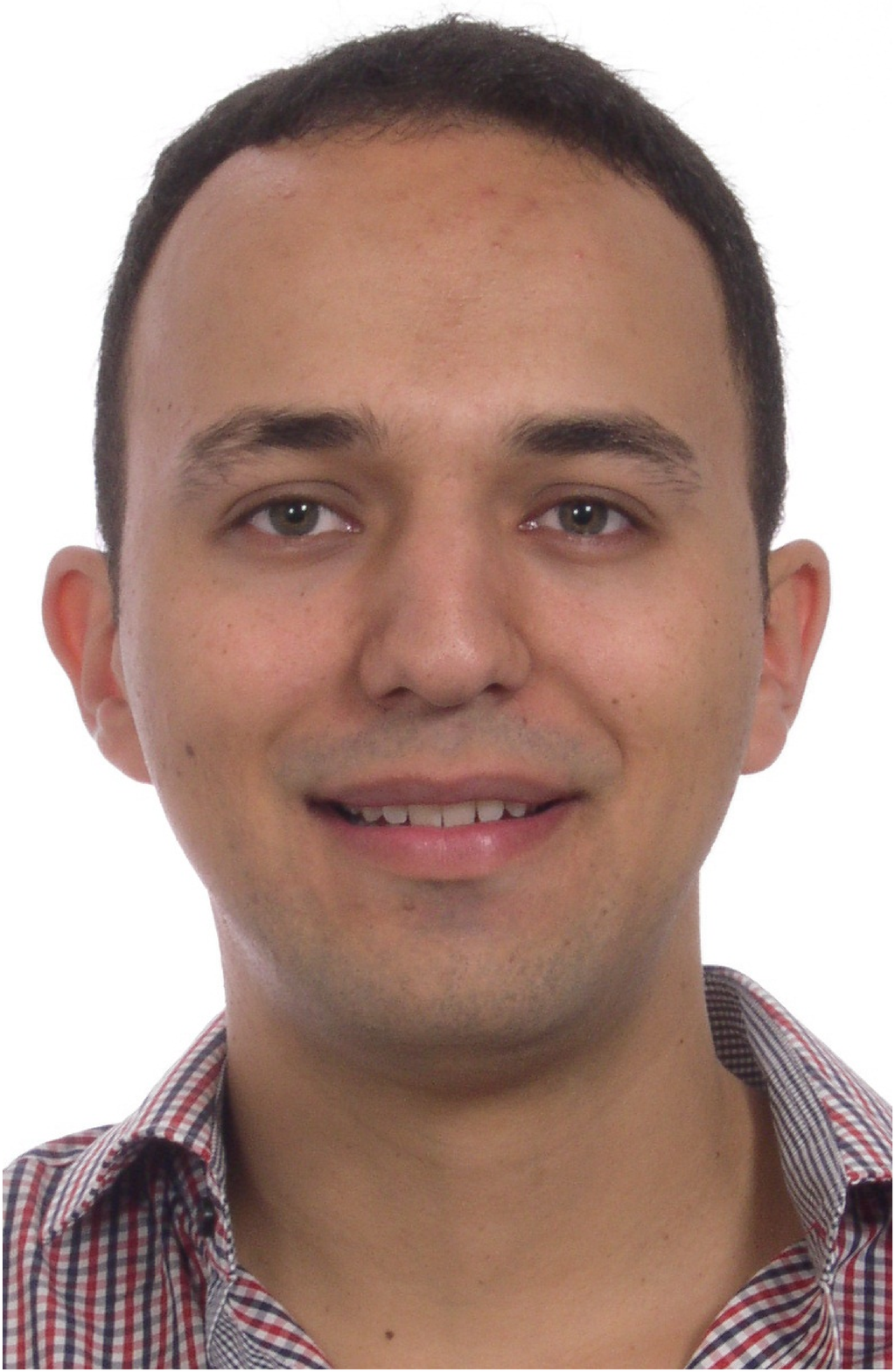}}]{Daniel Alberto Burbano Lombana}was born in San Juan de Pasto - Colombia. He received the B.S. degree in Electronic Engineering and the M.S. degree in Industrial automation from the National University of Colombia, in 2010 and 2012 respectively. He is currently a Ph.D. student in Computer and Automation Engineering at the University of Naples Federico II, Italy. His research is focused on distributed control of complex networks with applications to power grids.
\end{IEEEbiography}
\begin{IEEEbiography}[{\includegraphics[width=1in,height=1.5in,clip,keepaspectratio]{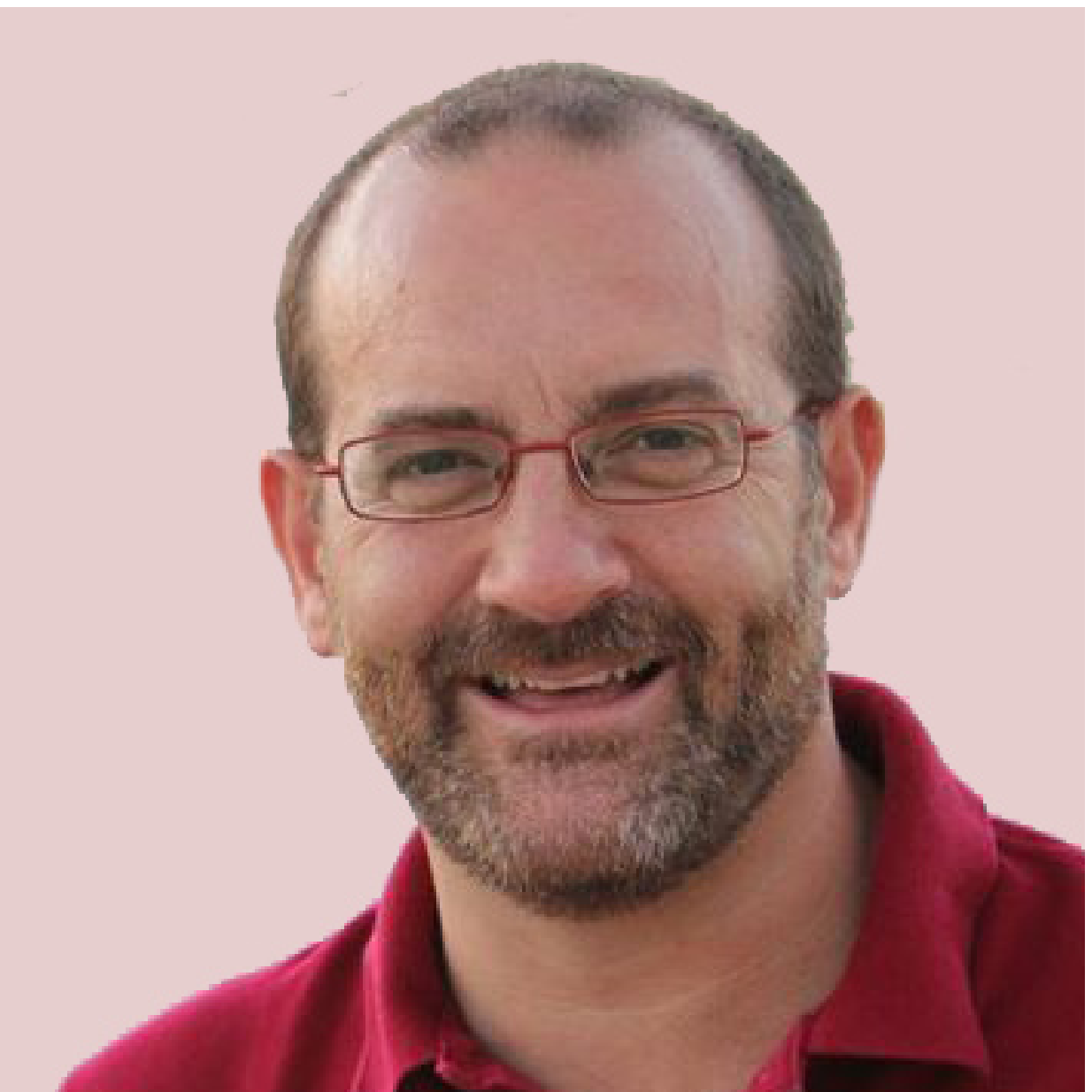}}]
{Mario di Bernardo} (SMIEEE 2006, FIEEE 2012) is Professor of Automatic Control at the University of Naples Federico II, Italy. He is also Professor of Nonlinear Systems and Control at the University of Bristol, U.K. In 1998, he obtained a Ph.D. in Nonlinear Dynamics and Control from the University of Bristol, U.K. He was then appointed to a Lecturership at the Department of Engineering Mathematics of the University of Bristol in 1997 where he became Professor of Nonlinear Systems and Control (part-time) on 1st August 2007. From 2001 till 2003, he was Assistant Professor at the University of Sannio, Italy.
On 28th February 2007 he was bestowed the title of Cavaliere of the Order of Merit of the Italian Republic for scientific merits from the President of Italy, HE Giorgio Napolitano. In January 2012 he was elevated to the grade of Fellow of the IEEE for his contributions to the analysis, control and applications of nonlinear systems and complex networks.
In 2006 and again in 2009 he was elected to the Board of Governors of the IEEE Circuits and Systems Society (one of the largest in the IEEE with over 9000 members). From 2011 he is serving as Vice President for Financial Activities of the IEEE Circuits and Systems Society. He authored or co-authored more than 200 international scientific publications including more than 100 papers in scientific journals, over 100 contributions to refereed conference proceedings, a unique research monograph on the dynamics and bifurcations of piecewise-smooth systems published by Springer-Verlag and two edited books. 
He serves on the Editorial Board of several international scientific journals and conferences. He is Deputy Editor-in-Chief of the IEEE Transactions on Circuits and Systems: Regular Papers. He is also Associate Editor of the IEEE Transactions on Control of Network Systems, Nonlinear Analysis: Hybrid Systems and Associate Editor of the Conference Editorial Board of the IEEE Control System Society and the European Control Association (EUCA).
He is regularly invited as Plenary Speakers in Italy and abroad and has been organizer and co-organizer of several scientific initiatives and events. He received funding from several agencies and industry.
\end{IEEEbiography}
\end{document}